\documentstyle[amscd,amssymb,verbatim,diagrams,12pt]{amsart}
\newcommand{\bC}{{\bf C}}

\newcommand{\vareps}{\varepsilon}

\newcommand{\ZZ}{{\cal Z}}

\renewcommand{\mod}{\operatorname{mod}}

\newcommand{\OO}{{\cal O}}

\newcommand{\ff}{\operatorname{for}_3}

\newcommand{\be}{{\bf e}}

\newcommand{\G}{{\Bbb G}}

\newcommand{\mg}{{\frak m}}

\newcommand{\UU}{{\cal U}}
\newcommand{\WW}{{\cal W}}

\newcommand{\Spec}{\operatorname{Spec}}

\newcommand{\Proj}{\operatorname{Proj}}
\renewcommand{\P}{{\Bbb P}}

\newcommand{\si}{\sigma}

\newcommand{\Pic}{\operatorname{Pic}}

\newcommand{\ga}{\gamma}
\newcommand{\de}{\delta}
\newcommand{\eps}{\epsilon}

\newcommand{\A}{{\Bbb A}}

\numberwithin{equation}{subsection}

\newtheorem{thm}{Theorem}[section]
\newtheorem{prop}[thm]{Proposition}
\newtheorem{lem}[thm]{Lemma}
\newtheorem{cor}[thm]{Corollary}
{  \theoremstyle{definition}

\newtheorem{rem}[thm]{Remark}

}

\newcommand{\Pf}{\noindent {\it Proof}}

\newcommand{\ov}{\overline}

\newcommand{\cusp}{\operatorname{cusp}}

\newcommand{\MM}{{\cal M}}

\newcommand{\Res}{\operatorname{Res}}

\newcommand{\Aut}{\operatorname{Aut}}

\renewcommand{\a}{\alpha}
\renewcommand{\b}{\beta}
\newcommand{\om}{\omega}

\newcommand{\la}{\lambda}

\newcommand{\C}{{\Bbb C}}

\newcommand{\R}{{\Bbb R}}
\newcommand{\Z}{{\Bbb Z}}
\newcommand{\Q}{{\Bbb Q}}

\newcommand{\wt}{\widetilde}

\newcommand{\sub}{\subset}
\newcommand{\ed}{\qed\vspace{3mm}}

\newcommand{\bv}{{\bf v}}
\newcommand{\ba}{{\bf a}}

\title{Birational models of $\MM_{2,2}$ arising as moduli of curves with nonspecial divisors}
\author{Drew Johnson}
\author{Alexander Polishchuk}
\address{University of Oregon}
\address{University of Oregon, National Research University Higher School of Economics, and
Korea Institute for Advanced Study}

\begin{document}

\begin{abstract}
We study birational projective models of $\MM_{2,2}$ obtained from the moduli space of curves with nonspecial
divisors. We describe geometrically which singular curves appear in these models and show that one of them is obtained
by blowing down the Weierstrass divisor in the moduli stack of $\ZZ$-stable curves $\ov{\MM}_{2,2}(\ZZ)$ defined by 
Smyth. As a corollary, we prove projectivity of the coarse moduli space $\ov{M}_{2,2}(\ZZ)$.
\end{abstract}

\maketitle

\section*{Introduction}
There has been a lot of interest recently in studying birational models of the moduli spaces of pointed curves $\MM_{g,n}$,
in particular, in connection with the Hassett-Keel program (see e.g., \cite{Smyth-modular}, \cite{FS}, \cite{AFSW} and references therein).
Typically, such birational models are constructed either by looking at a moduli functor defined
by some geometric restrictions on types of curves or as GIT quotients. In the former case projectivity of the coarse moduli
could be difficult to check, while in the latter case it is not clear how to describe geometrically
which curves are GIT stable. 

One example of a situation where it is possible to describe geometrically all GIT stable curves was considered
for certain moduli spaces of genus $1$ curves with $n$ marked points in \cite{P-g1}. 
In the present paper we study a similar picture for genus $2$ curves with $2$ marked points.
More precisely, we use the result of \cite{P-ainf} stating that there is an affine scheme $\wt{\UU}^{ns}_{2,2}$ with $\G_m^2$-action
parametrizing data $(C,p_1,p_2,v_1,v_2)$, where $C$ is a reduced connected curve of arithmetic genus $2$,
$p_1$ and $p_2$ are distinct smooth points of $C$, $v_1$ and $v_2$ are nonzero tangent vectors at them,
such that the line bundle $\OO_C(p_1+p_2)$ is ample and nonspecial (has vanishing $H^1$).

We get projective models of $\MM_{2,2}$ by considering GIT quotients of $\wt{\UU}^{ns}_{2,2}$ by $\G_m^2$.
The corresponding choices of stability were described in \cite{P-ainf}
(in a more general case of curves of genus $g$ with $g$ marked points): the space of characters of $\G_m^2$ contains
$5$ chambers for which the GIT quotient is nonempty. Up to permuting the marked points we have three 
stabilities to study, which we call (i),(ii) and (iii)-stabilities (where (i)-stability is invariant under the permutation of the marked
points). 

Our first result is an explicit description of these stabilities geometrically (see Theorem \ref{thm:geo-stab1}).
The description is case by case, depending on the types of irreducible components of a curve
and their intersection subscheme.
For (ii)- and (iii)-stabilities, we also give a short description using some cohomological invariants and the $e$-invariant
of singularities (see Theorem \ref{thm:geo-stab2}). 
In addition, we show that the moduli space for (iii)-stability is simply the weighted projective stack
$\P(1,2,3,3,4,5)$.

Our second main result (which is an analog of \cite[Thm.\ 2.4.5]{P-contr}) 
gives the relation between the moduli space for (i)-stability, $\ov{\UU}^{ns}_{2,2}$, and
the moduli of $\ZZ$-stable curves $\ov{\MM}_{2,2}(\ZZ)$ constructed by Smyth in \cite{Smyth-modular}
(here $\ZZ$ is the extremal assignment, assigning to a DM-stable curve its unmarked components).
Namely, we show that there is a regular morphism
$$\ov{\MM}_{2,2}(\ZZ)\to\ov{\UU}^{ns}_{2,2}$$
which blows down to a point the closure $\ov{\WW}$ of the Weierstrass divisor (consisting of $(C,p_1,p_2)$ such that $h^1(p_1+p_2)\neq 0$) and is an isomorphism away from it (see Theorem \ref{thm:Z}).

As a biproduct of our results we prove that Smyth's moduli space $\ov{\MM}_{2,2}(\ZZ)$ has a projective coarse moduli
(see Theorem \ref{projective-thm}).
In addition, we show that the map of forgetting $p_2$ induces an isomorphism of
the Weierstrass divisor $\ov{\WW}\sub\ov{\MM}_{2,2}(\ZZ)$ with the GIT quotient stack
$\ov{\UU}^{ns}_{2,1}(2)$ of genus $2$ curves $(C,p_1)$ such that $h^1(2p_1)=0$. The latter stack was shown to be isomorphic to
$\P(2,3,4,5,6)$ in \cite[Prop.\ 2.1.1]{P-contr} (see Proposition \ref{Weier-prop}).


\medskip

{\it Convention}. Everywhere in this paper we work over $\Spec(\Z[1/6])$. One exception is Theorem \ref{projective-thm}, where we work over $\Spec(\Q)$.

\medskip

{\it Acknowledgments}. 
A.P. is supported in part by the NSF grant DMS-1700642 and by the Russian Academic Excellence Project `5-100'. 
While working on this project, A.P. was visiting Korea Institute for Advanced Study and ETH Z\"urich. 
He would like to thank these institutions for hospitality and excellent working conditions.

\section{Moduli Spaces of Curves with Nonspecial Divisors}\label{gen-2-sec}

\subsection{Canonical generators and canonical parameters}

We recall some definitions and results from \cite{P-ainf} and \cite{P-krich}. 
For ${\mathbf a} = (a_1, \cdots, a_n)$, where $a_i\in \Z_{\ge 0}$ and $a_1 + \cdots + a_n = g$, we consider
the stack $\UU_{g,n}^{ns}({\mathbf a})$ parametrizing data $(C, p_1, \dots, p_n)$ where $C$ is a reduced (but possibly singular) genus $g$ curve with smooth points $p_1, \dots, p_n \in C$ such that the divisor $D = a_1p_1 + \cdots + a_np_n$ is nonspecial, that is, $H^1(C,D)=0$. We denote by $\wt{\UU}^{ns}_{g,n}({\mathbf a})$ the ${\mathbb G}_m^n$-torsor over $\UU_{g,n}^{ns}({\mathbf a})$ corresponding to choices of nonzero tangent vectors at the marked points. 
It is proved in \cite{P-ainf} that $\wt{\UU}^{ns}_{g,n}({\mathbf a})$ is an affine scheme of finite type over $\Spec(\Z[1/6])$. 

In the case  ${\mathbf a} = (1,\ldots,1)$ we simply write 
$\UU^{ns}_{g,g}:=\UU^{ns}_{g,g}(1,\ldots,1)$, $\wt{\UU}^{ns}_{g,g}:=\wt{\UU}^{ns}_{g,g}(1,\ldots,1)$.

We will make use of $(f_i[-m])$, the \emph{canonical generators}  for the ring $\OO(C - \{p_1,\ldots,p_n\})$, as well as $(t_{i,m})$, the \emph{canonical parameters of order $m$} at each of the marked points $p_i$.
 
For $a_i = 1$, these are constructed as follows. Select a parameter $t_{i,1}$ at $p_i$ that is compatible with the chosen tangent vector, that is $\left<v_i, t_{i,1}\right>=1$. Then, since $H^0(C,p_i + D)/H^0(C, \OO)$ has dimension 1, there is a unique function $f_i[-2]$ (up to an additive constant) with the expansion at $p_i$ of the form $t_{i,1}^{-2} + c_1t_{i,1}^{-1} + \cdots$. We now define $t_{i,2} = t_{i,1} - \frac{c_1}{2} t_{i,1}^2$ so that $f_i[-2]= t_{i,2}^{-2} + O(t_{i,2}^0)$.

Now, inductively, for $m \ge 3$, since $H^0(C,(mp_i + D)/H^0(C,((m-1)p_i + D)$ is one-dimensional, we have a unique function $f_i[-m]$ in the quotient whose expansion at $p_i$ begins with $t_{i,m}^{-m}$. One can lift this uniquely to $H^0(C,mp_i + D)/H^0(C,\OO)$ by requiring that the expansion has the form $t_{i,m-1}^{-m} + c_{m-1} t_{i,m-1}^{-1} + \cdots$. We then set 
$$t_{i,m} = t_{i,m-1} - \frac{c_{m-1}}{m} t_{i,m-1}^{m},$$ 
so that we have 
$$f_i[-m] = t_{i,m}^{-m} + O(t_{i,m}^0).$$

The formal limit of $t_{i,m}$ gives the canonical formal parameter $t_i$ at $p_i$.

Notice that in order to construct parameters of order $m$, we should work over a ring where every integer $\le m$ 
is invertible (and to define $t_i$ one should work over $\Q$). Note also that the functions $f_i[-m]$ are uniquely
defined up to additive constants.

 We give names to the coefficients of the expansions of the $f_i[-m]$ in the canonical parameters by writing
\begin{equation}
	f_i[-m] = t_i^{-m} + \sum_{k \ge 0} \alpha_{ii}[-m,k] t_i^k \label{eq:aij1}
\end{equation}
at $p_i$, and for $i \neq j$, we expand $f_i[-m]$ at $p_j$:
\begin{equation}
	f_i[-m] = \sum_{k \ge -1} \alpha_{ij}[-m,k]t_j^{k}. \label{eq:aij2}
\end{equation}
These coefficients are well defined except for $\alpha_{ii}[-m,0]$ and $\alpha_{ij}[-m,0]$, since there is an ambiguity in an additive constant when constructing the $f_i[-m]$.
One could resolve the ambiguity  by requiring, for example, that $\alpha_{ii}[-m,0]=0$.  We will make a slightly different choice for the case of $g=2$, $n=2$ in the next section, which is the only case where it will be relevant for us.

We will view the $\alpha_{ij}[-m,k]$ as functions on the moduli scheme $\wt{\UU}^{ns}_{g,n}({\mathbf a})$.

Notice that in order to compute $\alpha_{ij}[-m,k]$, one only needs canonical parameters of order $k+2$, and to get $\alpha_{ii}[-m,k]$, one needs order $m+k+1$.  

This construction can be easily generalized to the case $a_i \ge 2$. In the case that $a_i=0$, the function $f_i[-1]$ is still well defined up to an additive constant. However, there is some ambiguity in determining the canonical parameters at $p_i$. This can be resolved as in \cite{P-krich}. In this paper, however, when $a_i=0$, we will only need the function $f_i[-1]$ and not the canonical parameter.

Each moduli scheme $\ov{\UU}^{ns}_{g,n}(\ba)$ is equipped with a natural $\G_m^n$-action, rescaling the tangents vectors.
We denote by $[C^{\cusp}(\ba)]$ the unique $\G_m^n$-invariant point.
For example, $C^{\cusp}(1,1)$ is the union of two cuspidal curves of genus $1$, glued transversally at the cusp.
The curve $C^{\cusp}(2)$ is the one-point projective completion of $\Spec(k[t^3,t^4,t^5])$.


\subsection{Explicit description of $\wt{\UU}^{ns}_{2,2}$}\label{g-2-n-2-a-1-1-sec}

We are going to describe the scheme $\wt{\UU}^{ns}_{2,2}$ and the universal curve over it by explicit equations involving
some of the coordinates $\a_{ij}[p,q]$ (see \eqref{eq:aij1} and \eqref{eq:aij2}). 

Note that we can define all of the coordinates $\a_{ij}[p,q]$ as functions
on $\wt{\UU}^{ns}_{g,g}(1,\ldots,1)$ only if we work over $\Q$. 
However, we will see that in our case we can work over $\Z[1/6]$. 

The following result from \cite[Sections 1,2]{P-ainf} gives a nice criterion for canonical generators and determines a convenient way to normalize the constants for $f_i$ and $h_i$.

\begin{thm} \label{thm:fh}
	Let $(C, p_1, p_2, v_1, v_2)$ be a point of $\wt{\UU}^{ns}_{2,2}$ over some commutative ring $R$, such that
	$6$ is invertible in $R$.
	For $i=1,2$, let $t_i$ be formal parameters at $p_i$, compatible with $v_i$, and let 
	$f_i \in H^0(C,p_i + D)$ and $h_i \in H^0(C,2p_i + D)$, where $D=p_1+p_2$, be such that
	\[
	f_i = t_i^{-2} +O(t_i^{-1}), \;\;\; h_i = t_i^{-3} + O(t_i^{-2})
	\]
	Then the following are equivalent:
	\begin{enumerate}
		\item $f_i = f_i[-2]$ and $h_i = f_i[-3]$, with constant terms normalized by $\alpha_{ii}[-2,0]=0$ and $2\alpha_{ii}[-3,0] = 3\alpha_{ii}[-2,1]$.
		\item We have 
		\begin{equation*}
			h_i^2 - f_i^3 \in H^0(C,3D), \;\;\; f_i h_i^2- f_i^4 \in H^0(C,4D).
		\end{equation*}
		\item  The affine curve $C\setminus\{p_1,p_2\}$ has defining equations of the form
		\begin{equation}\label{g-2-n-2-curve-eq}
		\begin{array}{l}
		f_1f_2=\a_{21}h_1+\a_{12}h_2+\ga_{21}f_1+\ga_{12}f_2+a,\\
		f_ih_j=\a_{ij}f_j^2+\ga_{ij}h_j+\b_{ji}h_i+
		(\a_{ji}\b_{ij}+\a_{ij}\ga_{ji})f_j
		+(\vareps_{ji}-\a_{ij}\a_{ji}^2)f_i+b_{ij},\\
		h_i^2=f_i^3+\pi_if_i-\a_{ij}^3h_j+(\b_{ij}^2-3\a_{ij}^2\ga_{ij})f_j+s_i,\\
		h_1h_2=\b_{12}f_2^2+\b_{21}f_1^2+\vareps_{12}h_2+\vareps_{21}h_1+\psi_{12}f_2+\psi_{21}f_1+u,
		\end{array}
		\end{equation}
		where in the second and third equations $(i,j)=(1,2),(2,1)$, and $(\pi_i, a, b_{ij}, \psi_{ij}, s_i, u)$ are in $R$.
		
	\end{enumerate}
\end{thm}

Here the equivalence of (2) and (3) follows from \cite[Thm.\ 1.2.4]{P-ainf}, as well as from
formulas of \cite[Sec.\ 2.2]{P-ainf}. The equivalence with (1) is easy to check.

 If $f_i$ and $h_i$ satisfy any of the parts of the theorem, we will say that they are \emph{the} canonical generators.
 
 \begin{rem} \label{rem:make-can-gen}
Recall from \cite{P-ainf} a method for producing canonical generators from arbitrary nonconstant
$\wt{f}_i \in H^0(C,p_i+D)$ and $\wt{h}_i \in H^0(C,2p_i + D)$. One first rescales them if necessary to obtain expansions beginning with $t_i^{-2}$ and $t_i^{-3}$. Then one sets 
$f_i = \wt{f}_i + a_i$ and $h_i = \wt{h}_i + b_i\wt{f}_i + c_i$. We see that 
\[
h_i^2 - f_i^3 = \wt{h}_i^2 - \wt{f}_i^3 + 2b_i \wt{h}_i \wt{f}_i + (b_i^2 - 3a_i)\wt{f}_i^2 \mod H^0(C, 3D)
\]
One sees that $\wt{h}_i^2 -\wt{f}_i^3 \in H^0(C,2p_i + 3D)$ and also that $\wt{h}_i \wt{f}_i$ and $\wt{f}_i^2$ form a basis for $H^0(C,2p_i + 3D)/H^0(C, 3D)$. It follows that $a_i$ and $b_i$ can be chosen uniquely so that $h_i^2 - f_i^3 \in H^0(C, 3D)$.

Once $a_i$ and $b_i$ are fixed, we have
\[
f_i h_i^2- f_i^4 = \wt{f}_i \wt{h}_i^2 - \wt{f}_i^4 + 2 c_i \wt{f}_i \wt{h}_i \mod H^0(C, 4D).
\]
Since $\wt{f}_i \wt{h}_i^2 \in H^0(C, p_i + 4D)$ and $H^0(C, p_i + 4D)/H^0(C, 4D)$ is one-dimensional and
is spanned by $\wt{f}_i \wt{h}_i$, we can choose $c_i$ uniquely so that $f_i h_i^2- f_i^4 \in H^0(C, 4D)$.
\end{rem}

 We adopt the following more compact notation for the coordinates (following \cite[Sec.\ 2.2]{P-ainf}):
\begin{equation}\label{coordinate-notation}
	\a_{ij}=\a_{ij}[-2,-1], \  \ga_{ij}=\a_{ij}[-2,0], \
	 \b_{ij}=\a_{ij}[-3,-1], \
	\vareps_{ij}=\a_{ij}[-3,0]
\end{equation}

 In other words, if $t_i$ is the canonical parameter at $p_i$, we have expansions, for $i\neq j$:
\begin{align}
	f_i &= \alpha_{ij}t_j^{-1} + \gamma_{ij} + \cdots \nonumber \\ 
	h_i &= \beta_{ij}t_j^{-1} + \epsilon_{ij} + \cdots \label{eq:expansions} 
\end{align}

If one is only interested in $\a_{ij}$ and $\b_{ij}$, one only needs canonical parameters to order one (that is, the parameter is compatible with $v_i$). If one is interested in $\ga_{ij}, \vareps_{ij}$, one only needs the canonical parameters to order two---that is, it is sufficient to find parameters $t_i$ so that the expansion of $f_i$ at $p_i$ has vanishing $t_i^{-1}$ term. 


One can substitute these expansions into the equations of \eqref{g-2-n-2-curve-eq} and compare coefficients of various powers of $t_i$ to see that the names $\a_{ij}, \b_{ij}, \ga_{ij}, \epsilon_{ij}$ have been chosen consistently. 

In order to generate the ring of functions on $\wt{\UU}_{2,2}$, we will need (in addition to $\a_{ij}$, $\b_{ij}$, $\ga_{ij}$, $\vareps_{ij}$) the coefficient $\pi_i$ appearing in \eqref{g-2-n-2-curve-eq}. Notice that the conditions of Theorem \ref{thm:fh} imply that $h_i^2-f_i^3 \in H^0(C,2p_i+3p_j)$. Given a  parameter $t_i$ compatible with $v_i$, we see that $\pi_i$ is the coefficient of $t_i^{-2}$ in the expansion of $h_i^2-f_i^3$.


As usual, the Buchberger's algorithm gives the relations between the coefficients which are equivalent to the condition that the elements $(f_i^n, f_i^nh_i)$
form a basis of $H^0(C\setminus\{p_1,p_2\},\OO)$. These relations take better form after changing $\pi_i$ to
$$
\wt{\pi_i} = \pi_i + 3\a_{ij}^2\ga_{ji} + 3\a_{ij}\a_{ji}\b_{ij} + 3\ga_{ij}^2
$$
for $(i,j)=(1,2),(2,1)$.

\begin{prop}\label{g-2-n-2-a-1-1-prop} Let us work over $\Z[1/6]$.
The moduli scheme $\wt{\UU}^{ns}_{2,2}$ is isomorphic to the locus in the affine space $\A^{10}$
with the coordinates 
\begin{equation}\label{U-22-11-coordinates}
\a_{12},\a_{21},\b_{12},\b_{21},\ga_{12},\ga_{21},\vareps_{12},\vareps_{21},\wt{\pi}_1,\wt{\pi}_2,
\end{equation}
where the matrix
$$\left(\begin{matrix} \a_{12} & \a_{21} & \b_{12} & \b_{21} \\ 2\vareps_{12} & 2\vareps_{21} & \wt{\pi}_1 &\wt{\pi}_2 
\end{matrix}\right)$$
has rank $\le 1$ (i.e., the equations are given by the vanishing of all the $2\times 2$-minors). In other words,
$$\wt{\UU}^{ns}_{2,2}\simeq C(\P^1\times\P^3)\times \A^2,$$
where $C(\P^1\times\P^3)$ is the affine cone over $\P^1\times\P^3$ in the Segre embedding.
The weights of the $\G_m^2$-action are given by 
$$wt(\a_{ij})=2\be_i-\be_j, \ wt(\b_{ij})=3\be_i-\be_j, \ wt(\ga_{ij})=2\be_i, \ wt(\vareps_{ij})= 3\be_i, \  wt(\wt{\pi}_i)=4\be_i.$$

The universal affine curve $C\setminus\{p_1,p_2\}$ over  $\wt{\UU}^{ns}_{2,2}$ is given by \eqref{g-2-n-2-curve-eq} with
\begin{align*}
&a = - \alpha_{12}^2 \alpha_{21}^2 - \gamma_{12} \gamma_{21} + \beta_{12}\beta_{21}, \\
&b_{ij} = 
	\a_{ij} \a_{ji}^2\ga_{ij} -2\a_{ij}^2\a_{ji} \b_{ji} -2 \a_{ij}\ga_{ji}^2 + 2\a_{ji}\b_{ij}\ga_{ji} + \b_{ji}\vareps_{ij} - \ga_{ij} \vareps_{ji}, \\
&\psi_{ij} = 3\a_{ij}\a_{ji}\ga_{ij} + \ga_{ji}\b_{ij}, \\
&s_{i} = \vareps_{ij}^2 - 3\a_{ij}^3\vareps_{ji}  - \ga_{ij}\wt{\pi}_{i}  + 2\ga_{ij}^3 + 2\b_{ij}^2\ga_{ji}
 - 3\a_{ij}^2\ga_{ij}\ga_{ji} + 3\a_{ij}\a_{ji}\b_{ij}\ga_{ij} - \a_{ij}^2\b_{ij}\b_{ji}  +3\a_{ij}^4\a_{ji}^2,
\\
&u = \b_{21}\wt{\pi}_1 - \vareps_{12}\vareps_{21}  + 2\a_{12}^2\a_{21}\vareps_{21}
-2\a_{12}^3\a_{21}^3  - 2\ga_{21}^2\b_{12} - 2\ga_{12}^2\b_{21}
+ \a_{12}\a_{21} (3\ga_{12}\ga_{21} - 2\b_{12}\b_{21}).  
\end{align*}
The scheme $\wt{\UU}^{ns}_{2,2}$ is irreducible of dimension $7$, Cohen-Macauley and normal. 
\end{prop}

\Pf . The equations are obtained by applying the Buchberger's criterion to our equations
\eqref{g-2-n-2-curve-eq}, where we consider the degree lexicographical order on the monomials with
$$\deg f_i=2, \ \deg h_i=3, \ \ h_1>h_2>f_1>f_2$$
(see e.g, \cite[Thm.\ 15.8]{Eis-CA}).
More precisely, we get equations between the coefficients occurring in \eqref{g-2-n-2-curve-eq}
by expressing the monomials 
$$f_if_jf_k, \ f_if_jh_k, \ f_if_jh_j,\ f_ih_j^2,\  f_ih_jh_k, \ f_ih_ih_j,\ h_ih_jh_k, \ h_ih_j^2$$
in terms of the basis $(f_i^n,f_i^nh_i)$ in two different ways. The fact that the resulting affine scheme is
isomorphic to $\wt{\UU}^{ns}_{2,2}$ follows from \cite[Thm.\ 1.2.4]{P-ainf}.

The geometric properties of the scheme $\wt{\UU}^{ns}_{2,2}$ 
follow from the well known properties of the variety of matrices of rank 
$\le 1$.
\ed





\section{GIT stabilities on $\wt{\UU}^{ns}_{2,2}$}

\subsection{Descriptions of unstable loci using coordinates}
For a pair of linearly independent vectors $\bv_1,\bv_2\in \R^2$ let us denote by $\bC(\bv_1,\bv_2)$ 
(respectively, $\ov{\bC}(\bv_1,\bv_2)$) the open cone
$\R_{>0}\bv_1+\R_{>0}\bv_2$ (respectively, the closed cone $\R_{\ge0}\bv_1+\R_{\ge0}\bv_2$).
According to \cite[Sec.\ 2.4]{P-ainf}, there are $5$ chambers with nonempty quotients in the GIT picture for the action
of $\G_m^2$ on $\wt{\UU}^{ns}_{2,2}$: 
$$\bC_0(-\be_1+2\be_2,-\be_1+3\be_2), \ \bC_0(-\be_1+3\be_2,\be_2), \ \bC_0(\be_2,\be_1), \ 
\bC(\be_1,3\be_1-\be_2), \bC(3\be_1-\be_2,2\be_1-\be_2).$$
One of them, namely $\bC_0=\bC_0(\be_1,\be_2)$ is invariant with respect to swapping $p_1$ and $p_2$.
The remaining $4$ chambers consist of two pairs 
that are swapped by the transposition of $p_1$ and $p_2$. After using the transposition, the following Proposition gives a description of all the chambers.

For the rest of the section, we adopt the following notation. Since the divisor on $C$ is required to ample, there must be a marked point on each component of $C$. In the case that $C$ is reducible, we let $C_i$ be the component with the marked point $p_i$ for $i=1,2$, $g_i$ be the arithmetic genus of $C_i$, $\xi$ be the intersection subscheme, and $\ell(\xi)$ its length. We have $g_1 + g_2 + \ell(\xi) - 1 = 2$ (our curves are reduced, so there are no higher Tors and hence $C_1.C_2 = \ell(\xi)$). Also note that if $\ell(\xi)=1$, we must have $g_1=g_2=1$, otherwise the divisor would be special.

In the next Lemma we discuss some particular types of the curves that will play a role in distinguishing GIT stabilities.

\begin{lem}\label{cusp-lem} Let $k$ be an algebraically closed field,
$(C,p_1,p_2)$ the reducible curve in $\UU^{ns}_{2,2}(k)$
with $g_1=1$, $g_2=0$, $\xi$ supported at one point $q$,
such that $q$ is a cusp on $C_1$. We have $C_1\setminus \{p_1\}=\Spec k[x,y]/(y^2-x^3)$.
and the isomorphism type of $(C,p_1,p_2)$ is determined by the tangent vector $v=\la \partial_x+\mu \partial_y$ 
to $C_1$ at $q$ corresponding to the embedding $\xi\sub C_1$. More precisely, there are three possibilities.

\noindent
(i) If $\la\neq 0$, $\mu\neq 0$, then $C$ is isomorphic to the curve $C^{\cusp}_{1,1}$, whose affine part is
$$(y^2=x^3,z=0)\cup (z=y^2,x=y)\sub\A^3$$
(the entire curve is obtained by adding one point at infinity on each component).

\noindent
(ii) If $\la=0$, $\mu\neq 0$, then $C$ is isomorphic to the curve $C^{\cusp}_{0,1}$, whose affine part is
$$(y^2=x^3,z=0)\cup (z=y^2,x=0)\sub\A^3.$$

\noindent
(iii) If $\la\neq 0$, $\mu=0$, then $C$ is isomorphic to the curve  $C^{\cusp}_{1,0}$, whose affine part is
$$(y^2=x^3,z=0)\cup (z=x^2,y=0)\sub\A^3.$$

In the cases (i) and (ii), $C$ has a singularity at $q$ which is analytically equivalent to the plane singularity
$x(y^2-x^3)=0$. In the case (iii), one has $\dim T_qC=3$. 
\end{lem}

\Pf . The gluing of $C_1$ and $C_2$ along $\xi$ is
determined by a surjective homomorphism $\phi:k[[x,y]]/(y^2-x^3)\to k[u]/(u^2)$, which is determined by a pair
$(\la,\mu)\neq (0,0)$ such that $\phi(x)=\la u$, $\phi(y)=\mu u$. Rescaling $x$, $y$ and $u$, we can assume that
$(\la,\mu)$ is either $(0,1)$, $(1,1)$, or $(1,0)$.

The equivalence of singularities in $C^{\cusp}_{1,1}$ and $C^{\cusp}_{0,1}$ is achieved by a change of
the coordinate $t\mapsto t+t^2/2$ on the normalization of $C_1$.
\ed

\begin{prop}\label{unstable-locus-prop}
	 
(i) For $\chi$ in the interior of $\bC_0$ the $\chi$-unstable locus in $\wt{\UU}^{ns}_{2,2}$
is the union of two irreducible components $Z_1$ and $Z_2$, intersecting in a single point $C^{\cusp}(1,1)$, where
$Z_i$, for $i=1,2$, is given by the equations
$$\a_{ij}=\b_{ij}=\ga_{ij}=\vareps_{ij}=\wt{\pi}_i=0,$$
where $j=3-i$. 

The generic point of $Z_1$ corresponds to a reducible curve $C$ with components $C_1$ and $C_2$ smooth with genera $g_1 = 0$ and $g_2 = 1$, and with $C_1$ tangent to $C_2$ at a single point.


\noindent
(ii) For $\chi$ in $\bC(\be_1,3\be_1-\be_2)$, the $\chi$-unstable locus is the union of three irreducible closed subsets:
$$P_1:\a_{12}=\b_{12}=\a_{21}=\b_{21}=0,$$
$$P_2:\a_{12}=\b_{12}=\vareps_{12}=\wt{\pi}_1=0,$$ 
$$R:\ga_{12}=\vareps_{12}=\wt{\pi}_1=\a_{21}=\b_{21}=\ga_{21}=\vareps_{21}=\wt{\pi}_2=0.$$

The generic point of $P_1$ corresponds to the nodal union of two smooth elliptic curves (with one marked point on each).

The generic point of $P_2$ corresponds to the union of a reducible curve $C$ with components $C_1$ and $C_2$ smooth with genera $g_1 = 0$ and $g_2 = 1$, meeting transversely at two points.

The generic point of $R$ is the curve $C^{\cusp}_{1,1}$ (see Lemma \ref{cusp-lem}).

\noindent
(iii) For $\chi$ in $\bC(2\be_1-\be_2,3\be_1-\be_2)$, the $\chi$-unstable locus is the union of three irreducible components:
$$W=P_1\cup W' :\a_{12}=0,$$
where $W'$ is given by $\a_{12}=\vareps_{12}=0$,
and
$$S:\a_{21}=\b_{ij}=\ga_{ij}=\vareps_{ij}=\wt{\pi}_i=0.$$

The generic point of $W'$ corresponds to the case when $C$ is smooth but $p_1$ is a Weierstrass point, i.e., $h^1(2p_1)\neq 0$.

The generic point of $S$ is the curve $C^{\cusp}_{1,0}$.

\end{prop}

\Pf.
First, we check the assertions about the components of the unstable loci. For (i), 
the weights of all the coordinates \eqref{U-22-11-coordinates} on $\wt{\UU}^{ns}_{2,2}$ belong to one of the
two closed cones $\ov{\bC}(2\be_1-\be_2,\be_1)$ and $\ov{\bC}(\be_2,2\be_2-\be_1)$.
If for some point $p$ of $\wt{\UU}^{ns}_{2,2}$
there exist coordinates with weights in both these cones that do not vanish at $p$ then the point $p$
is $\chi$-semistable. Conversely, if all the coordinates that do not vanish at $p$ have weights only in one of the cones
then $p$ is $\chi$-unstable. 

For (ii), let us split the weights of all the coordinates into two closed cones
$\ov{\bC}(3\be_1-\be_2,2\be_1-\be_2)$ and $\ov{\bC}(\be_1,2\be_2-\be_1)$. By the same argument as in (i),
for $\chi\in \bC(\be_1,3\be_1-\be_2)$, a point $p$ is $\chi$-unstable precisely when all the coordinates belonging to one of these
two closed cones vanish at $p$. In one case of $\ov{\bC}(\be_1,2\be_2-\be_1)$ we get the locus $R$. In the case of $\ov{\bC}(3\be_1-\be_2,2\be_1-\be_2)$ we get the locus $P$ given by $\alpha_{12} = \beta_{12}=0$. 
Imposing the conditions on minors from Proposition \ref{g-2-n-2-a-1-1-prop}, we see that in fact 
$P=P_1\cup P_2$, as given.

For (iii), we split the weights into those lying in $\ov{\bC}(\be_1,2\be_2-\be_1)$ and those on the ray generated
by $2\be_1-\be_2$ to obtain $W$ and $S$.

Next, we verify the claims about the generic points of the components of the unstable locus.


Let's first look at the locus $P_2$, so we have $\alpha_{12}=\beta_{12}=\epsilon_{12}=\tilde \pi_1 = 0$. As long as $\alpha_{21}, \beta_{21} \neq 0$ we use the substitutions
\begin{align*}
	f_2 &= \alpha_{21}x + \gamma_{21} \\
	h_2 &= \beta_{21}y + \epsilon_{21} \\
	f_1 &= z + \gamma_{12}
\end{align*}
The first line of \eqref{g-2-n-2-curve-eq} becomes $h_1 = xz$, so now transforming all variables into $x,y,z$, we obtain
\begin{align*}
	\beta_{21}z(x-y) &= 0 \\
	z(x^2-3\gamma_{12} - z) &= 0 \\
	\alpha_{21}^{3} x^{3} - 3 \alpha_{21}^{3} \gamma_{12} x +& 3 \alpha_{21}^{2} \gamma_{21} x^{2} -  \alpha_{21}^{3} x z - 9 \alpha_{21}^{2} \gamma_{12} \gamma_{21} -  \beta_{21}^{2} y^{2} - 3 \alpha_{21}^{2} \gamma_{21} z + 3 \beta_{21}^{2} \gamma_{12} + \cdots\\
	\cdots \alpha_{21} \tilde \pi_2 x &- 2 \beta_{21} \epsilon_{21} y + \beta_{21}^{2} z = 0
\end{align*}
So we see that there is a component that is a cubic in the plane $z=0$, with the cubic being
\begin{equation} \label{eq:cubic}
\alpha_{21}^{3} x^{3} - 3 \alpha_{21}^{3} \gamma_{12} x + 3 \alpha_{21}^{2} \gamma_{21} x^{2} - 9 \alpha_{21}^{2} \gamma_{12} \gamma_{21} -  \beta_{21}^{2} y^{2} + 3 \beta_{21}^{2} \gamma_{12} + 
\alpha_{21} \tilde \pi_2 x - 2 \beta_{21} \epsilon_{21} y = 0
\end{equation}
which is smooth for generic parameters.

The other component is the parabola $x^2=z + 3\gamma_{12}$ in the plane $x=y$. To make sure this satisfies the last equation, one can substitute out $z$ and $y$ and obtain
\[
\alpha_{21} \tilde \pi_2 x - 2 \beta_{21} \epsilon_{21} x = 0
\]
which is true because of the conditions on the minors in Proposition \ref{g-2-n-2-a-1-1-prop}.

To find the intersection of these two components, let us set $y=x$, $z=0$. Then the parabola has solutions $x = \pm \sqrt{3\gamma_{12}}$. We have already verified above that this is a solution of the cubic. So when $\gamma_{12} \neq 0$, we have an intersection at two distinct points.

Finally, note that the coordinate $z$ has a pole at the point at infinity of the parabola. This corresponds to a pole of $f_1$, so we conclude that $g_1 = 0$. This verifies the description of the generic point of $P_2$.


On $Z_1$, we have additionally that $\gamma_{12}=0$, so the points of intersection coincide, as desired.

In order to continue, it will be convenient to analyze the generic points by interchanging $p_1$ and $p_2$. Let $R'$ be the locus defined by interchanging $1$ and $2$ in the definition of $R$, and similarly for $S$ and $S'$.

On $R'$, we have that all parameters except $\alpha_{21}$ and $\beta_{21}$ are 0. Now the cubic \eqref{eq:cubic} becomes 
\[
\alpha_{21}^{3} x^{3} -  \beta_{21}^{2} y^{2}
\]
so we see that it has become a cusp, with the cusp point being the point (0,0,0), which is the point of intersection with the parabola.

On $S'$, all the variables except $\alpha_{21}$ are 0, so our change of variables for $y$ is not valid. So we leave $h_2$ as is and obtain that the system is equivalent to:  
\begin{align*}
	h_2 z &= 0 \\
	z(z-x^2) &=0 \\
	\alpha_{21}^3x^3 - \alpha_{21}^3xz - h_2^2 &= 0
\end{align*}
So there are two components: the cusp $\alpha_{21}^3x^3 - h_2^2 = 0$ in the plane $z=0$, and the parabola $z=x^2$ in the plane $h_2=0$. Again they intersect in one point at the cusp, but this time the Zariski tangent space is $3$-dimensional.

Now let us do $P_1$. We make the linear change of variables
\begin{align*}
	f_1 &= z + \gamma_{12}\\
	f_2 &= x + \gamma_{21}\\
	h_1 &= w + \epsilon_{12}\\
	h_2 &= y + \epsilon_{21}
\end{align*}
and our system becomes
\begin{align*}
	xz &= 0 \\
	yz &= 0 \\
	3 \gamma_{12} z^{2} + z^{3} + \tilde \pi_1 z - 2 \epsilon_{12} w -  w^{2} &= 0 \\
	3 \gamma_{21} x^{2} + x^{3} + \tilde \pi_2 x - 2 \epsilon_{21} y -  y^{2} &= 0 \\
	xw &= 0\\
	yw &= 0
\end{align*}
Now we can easily see that there are two components: a cubic in the plane $z=w=0$ defined by the 4th equation above, and a cubic in the plane $x=y=0$ defined by the third equation above. They meet transversely at the origin (because their containing planes do). The cubics are generically smooth.



On $W'$ the generic curve is a smooth. Using Lemma 2.3.3 of \cite{P-ainf}, one can see that $p_1$ is a Weierstrass point.
\ed

\subsection{Some geometric properties of GIT quotient stacks}

Now let $k$ be an algebraically closed field of characteristic $\neq 2,3$.
In this section we denote by $\ov{\UU}_{2,2}(i)$, $\ov{\UU}_{2,2}(ii)$ and $\ov{\UU}_{2,2}(iii)$ the GIT quotient stacks over $k$ corresponding to the (i)-, (ii)-, and (iii)-stabilites, respectively 
(elsewhere we also denote $\ov{\UU}_{2,2}(i)$ as $\ov{\UU}^{ns}_{2,2}$). 
Note that all of them are irreducible proper DM-stacks of dimension $5$.
Furthermore, it is easy to see that they are toric.

\begin{cor}
(i) The stacks $\ov{\UU}_{2,2}(ii)$ and $\ov{\UU}_{2,2}(iii)$ are smooth, while $\ov{\UU}_{2,2}(i)$ has the unique singular point $q$ corresponding to the union of two nodal irreducible curves of arithmetic genus $1$, glued transversally at the node
(one has $\Aut(q)=(\Z/2)^2$).

There is an isomorphism of $\ov{\UU}_{2,2}(ii)$ with the quotient by $\G_m^2$ of the open subset in the 
$\A^7$ with the coordinates $\a_{12},\b_{12},\a_{21},\b_{21},\ga_{12},\ga_{21},x$, where either
$(\a_{12},\b_{12})\neq 0$ or $(\a_{21},\b_{21},\ga_{12},\ga_{21},x)\neq 0$. The $\G_m^2$-weights of the coordinates
are given by
$$wt(\a_{ij})=2\be_i-\be_j, \ wt(\b_{ij})=3\be_i-\be_j, \ wt(\ga_{ij})=2\be_i, \ wt(x)=\be_1+\be_2,$$
where $j=3-i$.

There is an isomorphism of $\ov{\UU}_{2,2}(iii)$ with the weighted projective stack ${\mathbb P}(1,2,3,3,4,5)$.

\noindent
(ii) One has
$$\Pic(\ov{\UU}_{2,2}(i))=\Pic(\ov{\UU}_{2,2}(ii))=\Z^2, \ \ \Pic(\ov{\UU}_{2,2}(iii))=\Z.$$
All of these groups are generated by the line bundles associated with characters of $\G_m^2$.
One has $A_4(\ov{U}_{2,2}(i))=\Z^3$. The additional generator is the class of the divisor $P_1$.
\end{cor}

\Pf .
(i) The statement about singularities follows easily from the fact that the singular locus of $\wt{\UU}_{2,2}^{ns}$ is given
by the equations 
$$\a_{ij}=\b_{ij}=\ga_{ij}=\vareps_{ij}=\wt{\pi}_i=0 \ \text{  for }i=1,2.$$

To get the description of $\ov{\UU}_{2,2}(ii)$ we observe that on the complement of $P_1$ we have a well defined
function $x$ given by one of the quotients $2\vareps_{12}/\a_{12}$, $2\vareps_{21}/\a_{21}$, $\wt{\pi}_1/\b_{12}$
and $\wt{\pi}_2/\b_{21}$, on appropriate open subsets. 

To verify the last assertion, we observe that  $\a_{12}\neq 0$ on the (iii)-semistable locus. Thus, using the $\G_m^2$-action
we can normalize $\a_{12}$ to be $1$, and replace the group $\G_m^2$ by the subgroup $\G_m=\{(\la,\la^2)\}$.
Next, due to the equations given by the $2 \times 2$ minors in Proposition \ref{g-2-n-2-a-1-1-prop}, we can eliminate the variables $\vareps_{21}, \wt{\pi}_1$, and $\wt{\pi}_2$. This also eliminates all the relations. We are left with the quotient of 
${\mathbb A}^6\setminus \{0\}$ with coordinates $\alpha_{21}$, $\beta_{21}$, $\beta_{12}$, $\vareps_{12}$, $\gamma_{12}$, $\gamma_{21}$, by $\G_m$. One checks that the weights of $(\lambda,\lambda^2)$ on these are 3, 5, 1, 3, 2, 4 respectively..



\noindent
(ii) The cases of $\ov{\UU}_{2,2}(ii)$ and $\ov{\UU}_{2,2}(iii)$ follow easily from their explicit descriptions as quotients.

In the case of $\ov{\UU}_{2,2}(i)$, we have to compute the $\G_m^2$-equivariant Chow group of codimension $1$
(resp., Picard group) of 
$\wt{\UU}^{ns}_{2,2}\setminus Z$, where $Z=Z_1\cup Z_2$.
Now we use the isomorphism $\wt{\UU}^{ns}_{2,2}\simeq C(\P^1\times\P^3)\times \A^2$.
(see Proposition \ref{g-2-n-2-a-1-1-prop}). Let $0\in C(\P^1\times\P^3)$ be the vertex.
The complement $C(\P^1\times\P^3)\setminus 0$ can be identified with the complement to the zero section
in the line bundle $\OO(-1,-1)$ over $\P^1\times \P^3$. This easily implies that 
$$\Pic(C(\P^1\times\P^3)\setminus 0)=A_4(C(\P^1\times\P^3\setminus 0)=\Z,$$ 
with the generator given as preimage of a point in $\P^1$.
Hence, 
$$\Pic^{\G_m^2}((C(\P^1\times\P^3)\setminus 0)\times \A^2)=
A_4^{\G_m^2}((C(\P^1\times\P^3\setminus 0)\times\A^2)=\Z^3,$$ 
with the subgroup $\Z^2$ coming from the characters of $\G_m^2$, and another generator given 
by the class of $P_1$.
It follows that
$$A_4^{\G_m^2}(\wt{\UU}^{ns}_{2,2}\setminus Z)=A_4^{\G_m^2}(\wt{\UU}^{ns}_{2,2})=
A_4^{\G_m^2}((C(\P^1\times\P^3)\setminus 0)\times\A^2)=\Z^3.$$

On the other hand, it is well known that no multiple of $P_1$ is locally principal near any point
of $0\times\A^2\sub C(\P^1\times\P^3)\times\A^2$. Since this locus has nonempty intersection with
the complement to $Z$ (namely, $0\times (\G_m^2)$), it follows
that the usual Picard group of $\wt{\UU}^{ns}_{2,2}\setminus Z$ is trivial. Hence,
all line bundles on $\ov{\UU}_{2,2}(i)$ come from characters of $\G_m^2$.
Since the only global invertible functions on $\wt{\UU}^{ns}_{2,2}\setminus Z$ are constant,
this induces an isomorphism of $\Z^2$ with the $\G_m^2$-equivariant Picard group.
\ed


\section{Geometric characterizations of stabilities}
In this section we provide some geometric criteria one can use to determine whether a curve in $\UU^{ns}_{2,2}$ is stable with respect to one of our GIT stabilities. 

\subsection{Invariants of singularities}

Here we discuss some invariants of a reduced curve singularity $(C,q)$, so $\OO=\OO_{C,q}$ denotes such a local ring,
$\ov{\OO}$ is its normalization, $r$ is the number of branches, $\de=\dim(\ov{\OO}/\OO)$.
The main invariant we are interested in is $e$, the dimension of the smoothing component of the semiuniversal base
(see \cite{Gr}). For a quasihomogeneous singularity, it can be computed using the formula (see \cite[Thm.\ 2.5(3)]{Gr})
$$e=2\de-r+t,$$
where $t=\dim(\om/\mg\om)$, where $\om$ is the dualizing module, $\mg\sub \OO$ is the maximal ideal.
Note that for Gorenstein singularities we have $t=1$.
To compute $t$ in general one uses the normalization $\pi:\wt{C}\to C$. Then 
$\om_C$ can be identified with the subsheaf of rational $1$-forms $\eta$ on $\wt{C}$ such that
$$\sum_{p\in\pi^{-1}(q)}\Res_p(\pi^*(f)\eta)=0 \ \text{ for any } f\in \OO_C.$$ 

To understand the invariant $t$ for transversal unions of singularities, the following result is helpful.

\begin{lem}\label{lem:Gor-omega}
(i) Let $(C,q)$ be a reduced curve germ, $\pi:\wt{C}\to C$ the normalization.
Define $\wt{\om}_C$ to be the subsheaf of rational $1$-forms $\eta$ on $\wt{C}$ such that
$$\sum_{p\in\pi^{-1}(q)}\Res_p(\pi^*(f)\eta)=0 \ \text{ for any } f\in \mg.$$ 
Then we have an exact sequence
\begin{equation}\label{om-ext-seq}
0\to \om_C\to \wt{\om}_C\to \OO_q\to 0.
\end{equation}
If $C$ is singular and Gorenstein then the induced sequence
$$0\to \om_C|_q\to \wt{\om}_C|_q\to k\to 0$$
is still exact.

\noindent
(ii) Now assume that $(C,q)$ is the transversal union of $(C_1,q)$ and $(C_2,q)$. 
Then we have exact sequences
\begin{equation}\label{om-union-seq1}
0\to \om_C\to \wt{\om}_{C_1}\oplus \wt{\om}_{C_2}\to \OO_q\to 0,
\end{equation}
\begin{equation}\label{om-union-seq2}
0\to \om_{C_2}\to \om_C\to \wt{\om}_{C_1}\to 0.
\end{equation}
\end{lem}

\Pf . (i) The morphism $\wt{\om}_C\to \OO_q$ is given by the sum of residues at the points of $\pi^{-1}(q)$.
By definition, its kernel is $\om_C$. The surjectivity is clear: we can choose a rational $1$-form $\eta$ with a pole
of order $1$ at one of the points of $\pi^{-1}(q)$ and regular elsewhere.

Now assume that $C$ is Gorenstein and singular. We want to prove that the map $\om_C|_q\to \wt{\om}_C|_q$
is injective. Since $\dim \om_C|_q=1$, we just need to check that this map is nonzero. Assume it is zero. Then we have
$\om_C\sub \mg\wt{\om}_C$. On the other hand, we always have $\mg\wt{\om}_C\sub\om_C$, so we get
$\om_C=\mg\wt{\om}_C$. It follows that $\wt{\om}_C|_q\simeq k$, so $\wt{\om}_C\simeq \OO_C$. Hence,
$\mg\simeq\om_C$. But $\om_C\simeq \OO_C$ since $C$ is Gorenstein, so we deduce that $\mg\simeq \OO_C$,
i.e., $C$ is smooth, which is a contradiction.

\noindent
(ii) Since the maximal ideal of $\OO_C$ is the sum of maximal ideals of $\OO_{C_1}$ and $\OO_{C_2}$,
we obtain
$$\wt{\om}_C=\wt{\om}_{C_1}\oplus \wt{\om}_{C_2},$$
which gives \eqref{om-union-seq1}. The second sequence follows from this, using the exact sequences
\eqref{om-ext-seq} for $C_1$ and $C_2$.
\ed

Here are some computations of the invariant $e$.

\begin{lem} \label{lem:e} (i) For the coordinate cross in $n$-space, one has $t=n-1$, $e=2n-3$. E.g., for the node, one has $e=1$.

\noindent
(ii) For the elliptic $n$-fold point, one has $t=1$, $e=n+1$. E.g., for the genus-$1$ cusp one has $e=2$, 
and for the tacnode one has $e=3$.
For the union of two smooth branches that are glued at a point $q$ along a subscheme of length $3$ one has $e=5$.

\noindent
(iii) For $C$, which is the transversal union of a line with a Gorenstein singularity $C'$, one has
$t(C)=2$. If in addition, $C'$ is quasihomogeneous then $e(C)=e(C')+2$.
For example, for the transversal union of the genus-$1$ cusp with a line, one has $e=4$.

\noindent
(iv) For $C$, which is the transversal union of two Gorenstein singularities $C_1$ and $C_2$, one has
$t(C)=3$. If in addition, $C_1$ and $C_2$ are quasihomogeneous then $e(C)=e(C_1)+e(C_2)+3$.
For example, for the transversal union of two genus-$1$ cusps, one has $e=7$.


\noindent
(v) For $C$, which is a union of the genus-$1$ cuspidal curve $C_1$ with a line, glued along a length $2$ subscheme supported at the cusp $q$, we have two possible singularities. For the plane singularity $(y^2-x^3)x=0$ (occurring in $C^{\cusp}_{1,1}$
and $C^{\cusp}_{0,1}$) one has $e=5$, while for the singularity in $C^{\cusp}_{1,0}$ one has $e=6$.
\end{lem}

\Pf .
(i) In this case $\wt{C}$ is the disjoint union of the components $C_i$ of $C$ and sections of $\om_C$ correspond to
collections of $1$-forms in $\om_{C_i}(q)$ with the sum of residues equal to $0$ (here $q\in C_i$ is the common point of
the intersection in $C$). This immediately gives $t=n-1$.

\noindent
(ii) These are Gorenstein singularities, so $t=1$, $e=2\de-r+1$.

\noindent
(iii) Let $C=L\cup C'$. By Lemma \ref{lem:Gor-omega}, $\dim \wt{\om}_{C'}|_q=2$. On the other hand, $\wt{\om}_L=\om_L(q)$,
so $\dim \wt{\om}_L|_q=1$.
Hence, from the exact sequence \eqref{om-union-seq1} we get
$$t(C)=\dim \om_C|_q\ge \dim \wt{\om}_L|_q+\dim \wt{\om}_{C'}|_q-1=2.$$
On the other hand, from the exact sequence \eqref{om-union-seq2} we get
$$t(C)=\dim \om_C|_q\le \dim \om_{C'}|_q+\dim \wt{\om}_L|_q=2,$$
so $t(C)=2$.

\noindent
(iv) By Lemma \ref{lem:Gor-omega}, $\dim \wt{\om}_{C_1}|_q=\dim \wt{\om}_{C_2}|_q=2$. Hence, from the exact sequence
\eqref{om-union-seq1} we get
$$t(C)\ge \dim \wt{\om}_{C_1}|_q+\dim \wt{\om}_{C_2}|_q-1=3.$$
On the other hand, from the exact sequence \eqref{om-union-seq2},
$$t(C)\le \dim \om_{C_1}|_q+\dim \wt{\om}_{C_2}|_q=3,$$
so $t(C)=3$.

\noindent
(v) The plane singularity $(y^2-x^3)x=0$ is quasihomogeneous and Gorenstein, and hence, has $t=2$, $e=5$. 

In the case $C=C^{\cusp}_{1,0}$,
the completion of $\OO_C$ is the subring in $k[[t]]\oplus k[[u]]$ linearly spanned by $(t^{\ge 3},0)$, $(0, u^{\ge 2})$,
$(t^2,u)$ and $(1,1)$. Thus, (the completion of) $\om_C$ consists of the rational $1$-forms
$(\sum_{n\ge -3}a_nt^ndt, \sum_{m\ge -2}b_mu^mdu)$ such that
$$a_{-3}+b_{-2}=0, \ \ a_{-1}+b_{-1}=0.$$
The subspace $\mg\om_C$ is spanned by pairs of regular $1$-forms, as well as $(t^{-1}dt,-u^{-1}du)$.
It follows that $t(C)=2$, $e(C)=6$.
\ed

Now let us consider curves of arithmetic genus $2$ that occur in the moduli space $\UU^{ns}_{2,2}$. Note
that such curves have at most two irreducible components. For $C=C_1\cup C_2$, where $C_1$ and $C_2$ are irreducible components, we denote by $\xi=C_1\cap C_2$ the intersection subscheme. Note that if $g_i$ is the arithmetic genus of $C_i$
then the exact sequence
$$0\to \OO_C\to \OO_{C_1}\oplus \OO_{C_2}\to \OO_{\xi}\to 0$$
shows that $2=g_1+g_2+\ell(\xi)-1$.

\begin{lem}\label{sing-types-lem}
	The following are the only types of singularities of curves occurring in $\UU^{ns}_{2,2}$ with $e \ge 5$ (up to swapping
	$C_1$ and $C_2$).
	\begin{enumerate}
		\item For irreducible curves, the cusp singularity $\hat{\OO}_{C,q}=k[[t^3,t^4,t^5]]$ with $e=5$.
		\item For $g_1=g_2=1$, $\ell(\xi)=1$, we have three such singularities: the intersection point is nodal on both curves, $e=5$; nodal on one curve and cuspidal on the other, $e=6$; cuspidal on both curves, $e=7$.
		\item For $g_1=1$, $g_2=0$, $\ell(\xi)=2$, if the intersection is at a single point and the genus 1 curve has a cusp at the intersection, then $e=5$ or $e=6$ (according to Lemma \ref{lem:e}(v)).
		\item For $g_1=g_2=0$, $\ell(\xi) = 3$, if the intersection is a single point, we have $e = 5$.
	\end{enumerate}
\end{lem}

\Pf. 
The proof is by considering cases and using Lemma \ref{lem:e}.
For (1), we use the classification of irreducible curves of genus 2, see e.g. \cite[Sec.\ 2.3]{P-contr}. The cusp
$k[[t^3,t^4,t^5]]$ is quasihomogeneous, so the computation of $e$ is straighforward. Note that the other irreducible genus $2$
cusp $y^2-x^5=0$ is Gorenstein and has $e=4$.
For (2) and (3), we note that an irreducible genus 0 curve is smooth and an irreducible genus 1 curve can only have a node or a simple cusp. For (4), in the case when the intersection is a single point, we have a planar singularity, so $e=2\de-r+1=5$.
\ed

\subsection{Characterizations of stabilities}

\begin{thm} \label{thm:geo-stab1} Let $(C,p_1,p_2)$ be a point of $\UU^{ns}_{2,2}$. We have the following
characterizations of GIT stabilities.
	
{\bf (i)-stability.}
If $C$ is irreducible, it is stable.

If $\ell(\xi)=1$, it is stable if and only if the intersection point is at most nodal on both curves. 

If $\ell(\xi)=2$, it is stable if there are two distinct points of intersection, but unstable if there is only 1.

If $\ell(\xi)=3$, it is stable.

{\bf (ii)-stability.}
If $C$ is irreducible, it is stable.

If $\ell(\xi)=1$, it is unstable.

If $\ell(\xi)=2$, it is unstable in the following two cases: (a)
$g_1 = 0$; (b) $g_1=1$, $\xi$ is supported at one point $q$, and $C_1$ has a cusp at $q$. 
Otherwise, it is stable.

If $\ell(\xi)=3$, it is stable.

{\bf(iii)-stability.}
If $C$ is irreducible, it is stable if and only if $h^1(2p_1)=0$.

If $\ell(\xi)=1$, it is unstable.

If $\ell(\xi)=2$, it is unstable in the following cases: 
(a) $g_1=0$; (b) $g_1=1$, $\xi$ is contained in the smooth locus of $C_1$ and the degree 2 
divisor on $C_1$ associated to $\xi$ is linearly equivalent to $2p_1$;
(c) $g_1=1$, $\xi$ is supported at one point $q$ which is a node on $C_1$, and the nonconstant function in
$H^0(C_1,2p_1)$ has a constant restriction to $\xi$;
(d) $C\simeq C^{\cusp}_{0,1}$; (e) $C\simeq C^{\cusp}_{1,0}$.
Otherwise, it is stable.

If $\ell(\xi)=3$, it is stable.
\end{thm}

The above criteria are useful in practice, but a somewhat more elegant statement can be provided for (ii)- and (iii)-stability as follows. 




\begin{thm} \label{thm:geo-stab2}
	We have the following characterizations of (ii)- and (iii)-stabilities for $(C,p_1,p_2)\in\UU^{ns}_{2,2}$.
	
	{\bf(ii)-stability.} $C$ is stable if and only if $h^1(C,3p_1) = 0$ and either $h^1(C,3p_2) = 0$ or all singularities have $e \le 4$.
	
	{\bf(iii)-stability.} $C$ is stable if and only if $h^1(C,2p_1)=0$  and all singularities have $e \le 5$.
\end{thm}

\Pf.
We will prove Theorem \ref{thm:geo-stab2}, assuming Theorem \ref{thm:geo-stab1}. The proof of Theorem \ref{thm:geo-stab1} will occupy Section \ref{sec:geo-stab-proof}.

First we make some observations. Let $(i,j)$ be either $(1,2)$ or $(2,1)$. If $\alpha_{ij} =0$, then $f_i \in H^0(C, 2p_i)$, which implies $h^1(C,2p_1)>0$. Conversely, if $h^1(C,2p_i) >0$, then we have a non-constant section in $H^0(C,2p_i)$. Since $h^0(C,2p_i+p_j) = 2$, we see then that $H^0(C,2p_i)=H^0(C,2p_i + p_j)$, so $\alpha_{ij}=0$. By a similar argument, one sees that $\alpha_{ij}=\beta_{ij}=0$ is equivalent to $h^1(C,3p_i) > 0$.

Also, note that if $h^1(C,3p_i) \ge 1$, then $C$ is reducible by \cite[Lem.\ 2.4.4]{P-contr}. 

Now, assume that $C$ is (ii)-stable. Recall that $P_1$ and $P_2$ are the two components of the locus $P$ where $\alpha_{12}=\beta_{12}=0$ (see the proof of Prop.\ \ref{unstable-locus-prop}). 
Since $C$ is not in $P_1$ or $P_2$, we have that $h^1(C,3p_1)=0$. Let us now assume that $h^1(C,3p_2) > 0$. Then $C$ is reducible.
The exact sequence
$$0\to \OO_C(3p_1)\to \OO_{C_1}(3p_1)\oplus \OO_{C_2}\to \OO_\xi\to 0$$
shows that 
$$2=h^0(C,3p_1) = h^0(C_1, 3p_1) - \ell(\xi) + 1 \ge 5 - g_1 - \ell(\xi).$$ 
Thus, the only possibilities for $g_1$ and $\ell(\xi)$ are $g_1 = 0$ , $\ell(\xi) = 3$; or $g_1 = 1$, $\ell(\xi) = 2$ (in both cases
$g_2=0$). 
The first case has $h^1(C,3p_2)=0$, and the second has singularities with $e \le 4$, except for the case when $\xi$ is supported
at one point $q$ and $C_1$ has a cusp at $q$, which is (ii)-unstable by Theorem \ref{thm:geo-stab1}.
 
Assume that $C$ is not (ii)-stable. If it is in $P_1$ or $P_2$, then $h^1(C,3p_1)>0$ and we are done, so assume $h^1(C,3p_1)=0$ and $C$ is in $R$.  Then $\alpha_{21}=\beta_{21}=0$ so $h^1(C,3p_2)>0$ and $C$ is reducible. As before, this means that the only possibilities are $g_1=g_2=0$, $\ell(\xi)= 3$ and $g_1=1$, $g_2=0$, $\ell(\xi)=2$. By Theorem \ref{thm:geo-stab1}, since $C$ is unstable, we deduce that $g_1=1$, $\xi$ is supported at one point $q$, and $C_1$ has a cusp at $q$, hence 
$e \ge 5$.

For (iii)-stability, we first consider irreducible curves. Since we know that in this case $e \le 5$ (see Lemma \ref{sing-types-lem}), Theorem \ref{thm:geo-stab1} and Theorem \ref{thm:geo-stab2} say the same thing.

Assume now that $C$ is (iii)-stable and reducible. It is not in $W$, so $h^1(C,2p_1)=0$. Arguing as before, we get
$$1=h^0(C, 2p_1) = h^0(C_1, 2p_1) - \ell(\xi) + 1 \ge 4-g_1 - \ell(\xi),$$ 
which leads to two possibilities: $g_1=g_2=0$, $\ell(\xi) = 3$ or $g_1=1$, $g_2=0$, $\ell(\xi)=2$. From the classification of singularities, we see that $e\le 5$ unless $C\simeq C^{\cusp}_{1,0}$. But in this case
$C$ is (iii)-unstable by Theorem \ref{thm:geo-stab1}.
 
 Now assume $C$ is (iii)-unstable, reducible, and $h^1(C,2p_1) = 0$, i.e., $h^0(C,2p_1)=1$.
As before, we have two cases: either $g_1=g_2=0$, $\ell(\xi) = 3$, or $g_1=1$, $g_2=0$, $\ell(\xi)=2$. 
The former case is not possible since $C$ is (iii)-unstable. In the latter case, according the Theorem \ref{thm:geo-stab1},
we have to consider the cases (b)--(e). It is easy to see that in the cases (b), (c) and (d) one has $h^0(C,2p_1)=2$, so
they cannot occur. In the remaining case (e) one has $C\simeq C^{\cusp}_{1,0}$ so $e=6$.
\ed

\subsection{Proof of Theorem \ref{thm:geo-stab1}} \label{sec:geo-stab-proof}
Note that the condition $h^1(p_1) = 1$ is equivalent to the condition $\a_{12}=0$, while the condition $h^1(3p_1) = 1$ is equivalent to $\a_{12}=\b_{12}=0$.

We also note that stability depends only on the curve and the points, not on the choice of tangent vectors. Hence, in what follows, we are free to select local parameters to first order in any way that is convenient.



\subsubsection{$C$ is irreducible}

If $C$ is irreducible,  then by Proposition \ref{unstable-locus-prop} it is (i)-stable and (ii)-stable, and it is (iii)-stable 
if and only if $h^1(2p_1)=0$.


\subsubsection{Formulas for $f_1$ and $h_1$ when $g_1=1$} \label{g1=1} 
In the cases where $g_1=1$, we write the affine part of $C_1$ in Weierstrauss normal form as $y^2=x^3+ax+b$, where the point at infinity is the marked point $p_1$. (Notice that a different but equivalent choice of Weierstrauss normal form corresponds to a different choice of tangent vector at $p_1$.)

One can check that the functions $x$ and $y$ have poles of order 2 and 3, respectively, at $p_1$.
We have $y^2-x^3 = ax+b \in H^0(C_1,2p_1)$, so assuming that $x$ and $y$ can be suitably extended to $C_2$, we have $f_1|_{C_1} = x$ and $h_1|_{C_2}=y$.  In this case, we see that $\pi_1 = a$. 

\subsubsection{Case: $g_1 = 1$ and $g_2=1$ }
In this case we have $\ell(\xi)=1$, that is, the intersection is transversal at a point $q$ (not necessarily a node).

Since $H^0(C_i,p_i)=\mathbb C$, we see that $f_i$ and $h_i$ are constant on $C_j$ for $i \neq j$ (with values $\gamma_{ij}$ and $\epsilon_{ij}$, respectively). Hence $\a_{ij}=\b_{ij}=0$, so it is never stable for (ii) and (iii).

We claim that it is stable for (i) if and only if $q=C_1\cap C_2$ is at most nodal on both $C_1$ and $C_2$.

Indeed, we claim that the condition that $q$ is a cusp on $C_1$ is precisely the condition that our curve is in $Z_1$.
First, assume we are in $Z_1$. Then, one can check that the cubic $h_1^2=f_1^3$ in the plane $h_2 = \epsilon_{21}, f_2 = \gamma_{21}$ is a solution of \eqref{g-2-n-2-curve-eq}. This cubic is $C_1$. We have $f_1(q)=\ga_{12}=0$, so $q$ is the cusp on $C_1$.

Next, assume that $C_1$ has a cusp at $q$, so in the notation of Section \ref{g1=1}, we have $a=b=0$ and $q = (0,0)$ in $(x,y)$ coordinates.  Writing a function $g$ on $C$ as  $g = (g|_{C_1}, g|_{C_2})$, we have $f_1 = (x,0)$ and $h_1 = (y,0)$.  Hence $\gamma_{12}=\epsilon_{12}=0$. Furthermore, we have $\wt{\pi}_1=\pi_1 = a = 0$, so we are in $Z_1$.


\subsubsection{Case: $g_1=1$, $g_2=0$, $\xi$ is supported at two distinct points.}
The curves are glued transversally at two points $q_1$ and $q_2$. 

We claim that any such curve is stable for (i) and (ii).
Our claim follows from the assertions $\ga_{21}\neq 0$ and either $\a_{12} \neq 0$ or $\b_{12} \neq 0$. 

To see that $\ga_{21}\neq 0$, we first note that $f_2|_{C_1}=\ga_{21}$ and $h_2|_{C_1}=\epsilon_{ij}$ are constant (because $H^0(C_1,p_1)=\mathbb C$).
Let $z$ be an affine coordinate for $C_2 \cong {\mathbb P}^1$ with $p_2$ being the point at infinity, and the intersection points being $\pm 1$. Then we let
\[
f_2|_{C_2}=z^2 - \frac23, \ \ h_2|_{C_2}=z(z^2-1).
\]
These satisfy $f(-1)=f(1)$ and $h(-1) = h(1)$ and so extend to functions on $C$ in $H^0(C,2p_2+p_1)$ and $H^0(C,3p_2+p_1)$, respectively.   One computes that
$h_2^2-f_2^3= \frac13 z^2 - \frac{8}{27} \in H^0(C, 2p_2) \subset H^0(C, 3p_2)$, so $f_2$ and $h_2$ are the canonical generators. Hence,
$$\ga_{21}=f_2(0)=-\frac{2}{3}\neq 0.$$
Furthermore, we see that $z^{-1}$ is a canonical parameter of order 2.



Finally, we need to prove that we cannot have $\a_{12}=\b_{12}=0$.
Indeed, if this were the case then both $f_1$ and $h_1$ would be constant on $C_2$.
Hence, each of our generators $(f_1,h_1,f_2,h_2)$
would take the same value on $q_1$ as on $q_2$, which contradicts the fact that they generate the ring $\OO(C\setminus \{p_1,p_2\})$.

To study (iii)-stability, we also need to look more closely at $f_1$ and $h_1$.  We see that the ring $\OO(C \setminus \{p_1,p_2\})$ is isomorphic to the subring of $k[x,y]/(y^2=x^3+ax+b) \oplus k[z]$ consisting of pairs $(f(x,y),g(z))$ so that $f(A,B) = g(-1)$ and $f(C,D) = g(1)$ (where $(A,B)$ and $(C,D)$ are the points of intersection). 

One easily checks that $f_1 = (x, \frac12(C-A)z + \frac12(C+A))$ and $h_1 = (y, \frac12(D-B) + \frac12(D+B)z)$ are the canonical generators, so we get $\alpha_{12} = \frac12(C-A)$. This is on $W$ if and only if $\alpha_{12}=0$, which means that $C=A$, which means that the two points of intersections are the fiber of the elliptic cover from $C_1$ to ${\mathbb P}^1$. Finally, our curve is not in $S$ since $\gamma_{21} \neq 0$.



\subsubsection{Case: $g_1=1$, $g_2=0$, $\ell(\xi)$=2 and $\xi$ is supported at a single point $q$} 
The restrictions $f_2|_{C_1}$ and $h_2|_{C_1}$ are constant as before. Let $z$ be an affine coordinate for $C_2={\mathbb P}^1$, with $p_2$ the point at $z=\infty$ and $\xi$ supported at $z=0$ and $\xi = \Spec k[z]/z^2$.

Now the functions $z^2$ and $z^3$ have poles of order 2 and 3, respectively, at $p_2$, and vanish on $\xi$. Hence we define $f_2$ to be $z^2$ on $C_2$ and 0 and $C_2$, and $h_2$ to be $z^3$ on $C_2$ and 0 on $C_1$, and see that these are the canonical generators. Hence all the coordinates with subscripts starting with $2$ are 0,  so we are on the locus $Z_2$, thus 
(i)-unstable.


To study (ii)- and (iii)-stabilities we need to look at $f_1$ and $h_1$. As before, let $y^2=x^3+ax+b$ be the equation for $C_1$, and let $Az+B$ and $Cz+D$ be the restrictions of $x$ and $y$, respectively, to $\xi$. Hence, the ring $\OO(C \setminus \{p_1,p_2\})$ is isomorphic to the subring of $k[x,y]/(y^2=x^3+ax+b) \oplus k[z]$ consisting of pairs $(f(x,y),g(z))$ such that $f(Az+B,Cz+D) = g(z) \mod z^2$.  The point $(B,D) \in C_1$ is the point of intersection with $C_2$. 


We see that the pairs $f_1=(x, Az+B)$ and $h_1=(y, Cz + D)$ are the canonical generators, and we saw that $z^{-1}$ is the canonical parameter at $p_2$ to order 3.  
 Hence we can conclude that $\alpha_{12}=A$, $\beta_{12}=C$, $\gamma_{12}=B$ and $\epsilon_{12}=D$, and $\pi_1 = a$ and $\tilde \pi_1 = a + 3 \gamma_{12}^2$. 

First, we check (ii)-stability. We cannot have both $A$ and $C$ equal to 0, or the restriction from the elliptic curve would not generate the intersection. Hence we are not on $P_1$ or $P_2$. If we are on $R$ then the point of intersection is at $(x,y)=(0,0)$ and $a=0$. This forces $b=0$, so we are in the case of the cusp $y^2=x^3$, meeting a rational curve at the cusp point. Otherwise, the curve is (ii)-stable.

Next, we check (iii)-stability. If it is in $S$ then $C_1$ is the cuspidal cubic meeting the rational curve $C_2$ as above and additionally has $0=\beta_{21}=C$. That is, the function $y$ from $C_1\setminus\{p_1\}$ restricts to a constant on $\xi$, hence,
$C$ is isomorphic to $C^{\cusp}_{1,0}$.

If the curve is in $W$, we have $h^0(C,2p_1)=2$, so there exists a nonzero function on $C_1$ with zero restriction to $\xi$, 
which leads to the cases (b), (c) and (d).

If the curve is not in $S$ or $W$ it is (iii)-stable.

\subsubsection{Case: $g_1=0$ and $g_2=1$}
All the computations in the previous two sections are valid with the $1$ and $2$ subscripts interchanged. 
Hence the curve is still (i)-stable if and only if the intersection is supported at two distinct points.

In either case we also get $\alpha_{12}=\beta_{12}=0$, so the curve is in either $P_1$ or $P_2$ and hence (ii)-unstable. 
It is also in $W$ and hence (iii)-unstable.

\subsubsection{Case: $g(C_1)=0$ and $g(C_2)=0$}
We can assume that $p_1=\infty$ and $p_2=\infty$.  There are three cases, but we will see that they all give stable curves. 
\begin{enumerate}
\item $\xi=q_1\cup q_2\cup q_3$. 
We may assume that the points $0$, $1$, and $\la \in C_1$ are glued transversely to $0$, $1$, and $\mu \in C_2$, respectively. We must have $\la \neq \mu$, otherwise $h^0(p_1+p_2)>1$. The algebra $\OO(C \setminus \{p_1,p_2\})$ is isomorphic to the subring of $k[x] \oplus k[y]$ consisting of pairs $(f(x),g(y))$ such that
\[
	f(0) = g(0); \;\;\; f(1)=g(1); \;\;\; f(\la)=g(\mu)
\]
We see that $(x^2 + (c-1)x, cy)$, with $c= \frac{\la^2-\la}{\mu-\la}$ is such a pair. Also, it has a pole of order 2 at $p_1$ and order 1 at $p_2$, so we deduce that it is $f_1$ (up to a constant), so $\alpha_{12} = c$. Since $\la \neq 0,1$, we have $c\neq 0$. Symmetrically, $\alpha_{21} \neq 0$, so such a curve is stable.

\item $\xi$ is the union of a length $2$ subscheme at $q_1$ with the simple point $q_2$.

We may assume that the curves are glued along $k[x]/x^2 \cong k[y]/y^2$ and also at the points $\lambda \in C_1$ and $\mu \in C_2$, where $\lambda, \mu \neq \infty, 0$.

Then the ring $\OO(C \setminus \{p_1,p_2\})$ is isomorphic to the subring of $k[x] \oplus k[y]$ consisting of pairs $(f(x),g(y))$ satisfying $f(\lambda)=g(\mu)$ and $f(x) = g(x) \mod x^2$. Let $a= \frac{\lambda^2}{\mu-\lambda}$. The pair $(x^2 +ax, ay)$ satisfies this. This pair has a pole of order 2 at $p_1$ and order 1 at $p_2$, so this is, up to an additive constant, $f_1$. Hence $\alpha_{12} =a \neq 0$. Similarly, $\alpha_{21} \neq 0$, so these curves are always stable.

\item $\xi$ is a length $3$ subscheme at $q$ (we refer to this case as {\it the union of two osculating $\P^1$'s}). 
In this case we can assume that $\xi = \Spec k[x]/(x-\lambda)^3 \subset C_1$. Let $Ax^2 + Bx + C$ be the restriction of $y$ to  $k[x]/(x-\lambda)^3$, so the ring $\OO(C-\{p_1,p_2\})$ is isomorphic to
the subring in $k[x]\oplus k[y]$ consisting of pairs $(f(x),g(y))$ such that
$$f(x)\equiv g(Ax^2+ Bx + C) \mod (x-\lambda)^3,$$
 The condition that $h^0(p_1+p_2)=1$ implies that $A\neq 0$ (otherwise $(Bx + C, y)$ would be a non-constant section). Then we see that
that
\[
f_1 = (x^2 + \frac BA x + \frac CA,\frac 1A y)
\]
up to an additive constant,
which gives $\a_{21}=\frac 1A \neq 0$. Similarly, $\a_{12}\neq 0$, so these curves is always stable. 
\end{enumerate}

\section{Connection to $\ZZ$-stability}

In this section we only consider (i)-stability and refer to it simply as {\it GIT stability}.
We denote the corresponding GIT moduli stack by $\ov{\UU}^{ns}_{2,2}=\ov{\UU}^{ns}_{2,2}(1,1)$.

\subsection{The rational map $\ff$}

In this subsection we will use the abbreviation $\wt{\UU}(a_1,\ldots,a_n):=\wt{\UU}_{2,n}(a_1,\ldots,a_n)$ (where $n$ is either $2$ or $3$).

By $\wt{\UU}((1,0,1),(1,1,0))$ we denote the open subset of $\wt{\UU}(1,0,1)$ of curves $(C, p_1, p_2, p_3)$ that additionally satisfy $h^1(C, p_1 + p_2)=0$. Similarly, the subsets  $\wt{\UU}((1,0,1),(2,1,0))$ and $\wt{\UU}((1,0,1),(1,2,0))$ are defined by conditions $h^1(C, 2p_1 + p_2)=0$ and $h^1(C, p_1 + 2p_2)=0$, respectively.

There is a forgetful map
\begin{equation}\label{ff-mor-eq}
\ff: \wt{\mathcal U}((1,0,1),(1,1,0)) \rightarrow \wt{\mathcal U}(1,1)
\end{equation}
obtained by forgetting the point $p_3$ and the tangent vector at it.

Let $Z \subset \wt{\mathcal U}_{2,3}((1,0,1),(1,1,0))$ be the inverse image under $\ff$ of the unstable locus (with respect to 
(i)-stability), giving us
\[
\ff': \wt{\mathcal U}_{2,3}((1,0,1),(1,1,0)) \setminus Z \rightarrow \overline{\mathcal U}_{2,2}(1,1)
\]

 Let $f_i[-p]$ be the canonical generators over $\wt{\UU}(1,0,1)$, as introduced in Section \ref{gen-2-sec}. For brevity of notation (and following \cite{P-contr}) we let 
 \[
\alpha = \alpha_{23}[-1,-1], \;\; \beta_i = \alpha_{i3}[-2,-1]  
\]
(for $i=1,2$).

\begin{prop}
	The open subset $\wt{\UU}((1,0,1),(1,1,0)) \subset \wt{\UU}(1,0,1)$ is given by $\alpha \neq 0$. 
	
	The subset $\wt{\UU}((1,0,1),(2,1,0)) \subset \wt{\UU}(1,0,1)$ is given by $\alpha \neq 0$ or $\beta_{1} \neq 0$.
	
	The subset $\wt{\UU}((1,0,1),(1,2,0)) \subset \wt{\UU}(1,0,1)$ is given by $\alpha \neq 0$ or $\beta_{2} \neq 0$.
	
	Hence the intersection ${\mathcal Y} := \wt{\UU}((1,0,1),(2,1,0),(1,2,0)) \subset \wt{\UU}(1,0,1)$ is given by $\alpha \neq 0$ or $\beta_1\beta_2 \neq 0$.
\end{prop}


\Pf . For any $C$ in $\wt{\UU}(1,0,1)$, we have $h^0(C,p_1+p_2+p_3)=2$, hence $H^0(C,p_1+p_2+p_3)$ is spanned by $1$ and $f_2[-1]$. We now have $h^0(C,p_1 + p_2) = 2$ if and only if $H^0(C,p_1 + p_2)$ is also spanned by $1$ and $f_2[-1]$, which is the case if and only if  $f_2[-1]$ is regular at $p_3$, which is precisely the condition that $\alpha=0$.

Similarly, we see that $H^0(C, 2p_1 + p_2 + p_3)$ is $3$-dimensional and spanned by $1, f_2[-1]$, and $f_1[-2]$. We will have $h^0(C,2p_1 + p_2)=3$ if and only if $f_2[-1]$ and $f_1[-2]$ are regular at $p_3$, which is precisely the condition that $\alpha = 0$ and $\beta_1=0$.

Finally, the space $H^0(C,p_1 + 2p_2 + p_3)$ has $1,f_2[-1],f_2[-2]$ as a basis and $h^0(p_1+2p_2)=3$ if and only if
$f_2[-1]$ and $f_2[-2]$ are regular at $p_3$, i.e., $\a=0$ and $\b_2=0$.
\ed

Let $\wt{\mathcal W}$ be the locus in $\mathcal Y$ where $h^1(p_1 + p_2)>0$, so $\wt{\mathcal W}$ is defined by $\alpha = 0$. 
Note that $\wt{\mathcal U}((1,0,1),(1,1,0))$ is an open subset of $\mathcal Y$ whose complement is $\wt{\mathcal W}$, and $Z \cap \wt{\mathcal W} = \emptyset$.

The following is the key step in proving that there is a regular morphism from the moduli of $\ZZ$-stable curves
to $\ov{\UU}_{2,2}$ (see Theorem \ref{thm:Z} below).

\begin{prop} \label{prop:forgetful}
There exists a regular morphism 
\[
 \wt{\ff}: {\mathcal Y}  \rightarrow \wt{\UU}(1,1)
\]
which agrees with $\ff$ (see \eqref{ff-mor-eq}) after passing to quotients by $\G_m^2$.
The image of $\wt{\mathcal W}$ under $\wt{\ff}$ avoids the unstable locus, and by composition with the quotient we get a map
\[
 \wt{\ff}': {\mathcal Y} \setminus Z \rightarrow \overline{\UU}_{2,2}(1,1)
\]
which extends $\ff'$. 

The map $\wt{\operatorname{for}}_3'$ sends $\wt{\mathcal W}$ to a single point, whose explicit coordinates are given. 
\end{prop}

\Pf. Let  $\alpha_{ij}[-m,-k]$ be the coefficients of the expansions of the canonical generators $f_i[-m]$
 in the canonical parameter $t_j$ at $p_j$ over $\wt{\UU}(1,0,1)$ (see Section \ref{gen-2-sec}).
 Under the map 
$\ff:\wt{\UU}((1,0,1),(1,1,0))\to\wt{\UU}_{2,2}(1,1)$, the coordinates $\a_{ij}, \b_{ij}, \ga_{ij}, \eps_{ij}, \pi_i$ can be written as
 rational functions in the $\alpha_{ij}[-m,k]$. We will show that a modified version of this map can be extended to the locus where $\alpha = 0$. We won't need to use coordinates of the form $\alpha_{ij}[-m,0]$, so we don't need to be concerned about normalizing the constants.


We are going to apply the standard procedure for computing the canonical generators $f_i, h_i$, $i=1,2$,
associated with the family of $2$-pointed curves over $\wt{\UU}((1,0,1),(1,1,0))$ given by the map $\ff$ (see
Remark \ref{rem:make-can-gen}).
To begin with we set, for $i=1,2$:
\begin{align*}
	\wt{f_i} &:= f_i[-2] - \frac{\beta_i}{\alpha} f_2[-1] \\\
	\wt{h_i} &:= f_i[-3] - \frac{\alpha_{13}[-3,-1]}{\alpha} f_2[-1]
\end{align*}
The coefficients of $f_2[-1]$ are chosen to cancel out the pole at $p_3$, so for $(i,j)=(1,2)$ or $(2,1)$ we have $\wt{f_i} \in H^0(2p_i + p_j)$ and $\wt{h_i} \in H^0(3p_i + p_j)$.

Let $\phi = \alpha_{21}[-1,-1]$, $\lambda_1 = \beta_1 \phi /\alpha$ and $\lambda_2 = \beta_2/\alpha$. 

Let us define an increasing filtration $F_n$ on the space of Laurent series in a variable $z$ with coefficients being regular functions on $\wt{\UU}((1,0,1),(1,1,0))$. By definition, a Laurent series is in $F_n$ if it can be written as
$\sum_i a_i\a^{-i-n}z^i$, where $a_i$ extends to a regular function on $\wt{\UU}(1,0,1)$.
For example, $z,\alpha\in F_{-1}$ and $\lambda_1,\lambda_2\in F_1$. 
We will use the same notation $F_n$ for Laurent series in any variable.
 
Using these definitions, we see that the expansions in the canonical parameters $t_i$ of $\wt{\UU}(1,0,1)$ are of the form:
\[
\wt{f}_i = t_i^{-2} - \lambda_i t_i^{-1} + {F_1}
\]
and 
\[
\wt{h}_i = t_i^{-3} + F_{2} 
\]
Now, as in Remark \ref{rem:make-can-gen}, we wish to find functions (on $\wt{\UU}((1,0,1),(1,1,0))$) $a_i,b_i,c_i$ so that the canonical generators on $\wt{\UU}(1,1)$ are given by
\[
f_i := \wt{f}_i + a_i, \;\; h_i := \wt{h}_i + b_i \wt{f}_i +c_i
\]
 We find that $a_i = \frac 34 \lambda_i^2 + F_1$, $b_i = -\frac32 \lambda_i+F_0$, and $c_i = -\frac 12 \lambda^3 + F_2$. 
Hence we have
\begin{align}
	f_i &= t_i^{-2} - \lambda_i t_i^{-1} + \frac34 \lambda_i^2 + F_1 \nonumber \\
	h_i &= t_i^{-3} - \frac32\lambda_it_i^{-2} + \frac32\lambda_i^2t_i^{-1} - \frac12 \lambda_i^3 + F_2 \nonumber \\
	h_i^2-f_i^3 &= -\frac{3}{16}\lambda_i^4t_i^{-2} + \cdots \label{eq:h2f3pi}
\end{align}
so we see that $\pi_i = -\frac{3}{16}\lambda_i^4$.

Next, we see that the substitution $t_i = u_i - \frac12 \lambda_i u_i^2 + F_{-2}$  gives us
\begin{equation} \label{eq:fihi}
	f_i = u_i^{-2} + O(u_i^0) +F_1
\end{equation}
so $u_i$ is a canonical parameter of order $2$.

Next, we compute the expansion of $f_i$ and $h_i$ at $p_j$ for $i \neq j$. First we note that 
\begin{align*}
\wt{f}_1 &= -\frac{\beta_1}{\alpha}t_2^{-1} + F_1 \\
\wt{f}_2 &= -\lambda_2 \phi t_1^{-1} + F_1 
\end{align*}
and the expansion of $\wt{h}_i$ at $p_j$ is in $F_2$. Hence we have
\begin{align*}
 f_1 &= \wt{f}_1 + a_1 =   -\frac{\beta_1}{\alpha} t_2^{-1} + \frac34\lambda_1^2 +  F_{1} \\
 &= - \frac{\beta_1}{\alpha}(u_2^{-1} + \frac12\lambda_2 + \cdots) + \frac34\lambda_1^2+ F_1\\
 f_2 &= \wt{f}_2 +a_2  = -\lambda_2 \phi t_1^{-1} +\frac34 \lambda_2^2  + F_{1}\\
 &= -\lambda_2 \phi(u_1^{-1} + \frac12\lambda_1 + \cdots) +\frac34 \lambda_2^2 + F_1\\
 h_1 &= \wt{h}_1 +b_1 \wt{f}_1 +c_1 = \frac32\lambda_1 \frac{\b_1}{\a}t_2^{-1} - \frac12 \lambda_1^3 + F_2 \\
 &=  \frac32\lambda_1 \frac{\b_1}{\a}(u_2^{-1} + \frac12\lambda_2 +  \cdots) - \frac12 \lambda_1^3 + F_2\\ 
 h_2 &= \wt{h}_2 + b_2 \wt{f}_2 +c_2 = \frac32\lambda_2^2\phi  t_1^{-1} - \frac12 \lambda_2^3 + F_2 \\
 &= \frac32  \lambda_2^2 \phi(u_1^{-1} + \frac12\lambda_1 + \cdots) - \frac12 \lambda_2^3 +  F_2
\end{align*}

Now one can read off the functions $\a_{ij}, \b_{ij}, \ga_{ij}, \eps_{ij}$ (resp., $\pi_i$)
from the above equations (resp., from \eqref{eq:h2f3pi}). 
We next wish to consider the modified map  $\wt{\ff}:=(\alpha, \alpha) \cdot \ff$ 
(where we use the action of $\G_m^2$ on $\wt{\UU}(1,1)$) and show that it extends to the locus where $\alpha = 0 $. To do this, we just need to check that the expressions for $\a_{ij}, \b_{ij}, \ga_{ij}, \eps_{ij}, \pi_i$ after the action do not have poles at $\alpha = 0$. We also wish to compute explicitly the coordinates of points in the image of $\wt{\WW}$. 
A consequence of our definition of $F_n$ is that after the action $(\alpha, \alpha)$, the terms of $f_i$ in $F_1$ and the terms of $h_i$ in $F_2$ will vanish along $\alpha = 0$.

Hence we see that after acting by $(\alpha, \alpha)$, the image of a point in $\wt{\WW}$ will have  
\begin{align*}
\alpha_{12} = -\beta_1, \;\;\; &\alpha_{21} = -\beta_2 \phi,\;\; \\
\beta_{12} = \frac32 \beta_1^2 \phi, \;\;\;\; &\beta_{21} = \frac32 \beta_2^2\phi \\
\gamma_{12} = -\frac12 \beta_1 \beta_2 + \frac34 \beta_1^2 \phi^2, \;\;\;\; &\gamma_{21} = -\frac12 \beta_1 \beta_2 \phi^2 + \frac34 \beta_2^2\\
\epsilon_{12} = \frac34 \phi \beta_1^2\beta_2 - \frac12\beta_1^3\phi^3, \;\;\;\; &\epsilon_{21} = \frac34 \beta_2^2\beta_1 \phi^2 - \frac12 \beta_2^3
\end{align*}

Also, from \eqref{eq:h2f3pi} we obtain
\[
\pi_1 =  -\frac{3}{16} \beta_1^4\phi^4,\;\;\;\;\; \pi_2 = -\frac{3}{16}  \beta_2^4.
\]

We can now see that these functions are regular on $\alpha=0$. Hence $\wt{\operatorname{for}}_3$ is well defined on $\mathcal Y$.

We need to check that no point of $\wt{\mathcal W}$ maps to the unstable locus in $\wt{\UU}_{2,2}(1,1)$. 

\begin{lem}
	On $\wt{\mathcal W}$, we have $\beta_2 = -\phi^2 \beta_1$.
\end{lem}
\Pf. 
We know that the space $H^0(C,2p_1 + 2p_2 + p_3)$ is $4$-dimensional, and $1, f_2[-1],f_2[-2],f_1[-2]$ form its basis. Since $f_2[-1]^2 \in H^0(2p_1 + 2p_2) \subset H^0(2p_1 + 2p_2 + p_3)$, we see that
\[
f_2[-1]^2 = a + bf_2[-1] + cf_2[-2] + df_1[-2].
\]
for some constants $a,b,c,d$.

Looking at the coefficient of $t_1^{-2}$ shows that $d = \phi^2$. Looking at the coefficient of $t_2^{-2}$ shows that $c=1$. Finally, since $\a=0$ over $\wt{\mathcal W}$, looking at the coefficient of $t_3^{-1}$ gives $0 = c\beta_2 + d\beta_1$ and the result follows. \ed

Now we can describe the image of a point of $\wt{\mathcal W}$ under the map $\wt{\operatorname{for}}_3$ as:
\begin{align*}
	\alpha_{12} = -\beta_1, \;\;\; &\alpha_{21} = \beta_1 \phi^3,\;\; \\
	\beta_{12} =\frac32\beta_1^2 \phi, \;\;\;\; &\beta_{21} = \frac32\beta_1^2\phi^5 \\
	\gamma_{12} = \frac54 \beta_1^2 \phi^2, \;\;\;\; &\gamma_{21} = \frac54 \beta_1^2 \phi^4\\
	\epsilon_{12} = -\frac54 \beta_1^3\phi^3, \;\;\;\; &\epsilon_{21} = \frac54 \beta_1^3 \phi^6\\
	\pi_1 = -\frac{3}{16} \beta_1^4 \phi^4 \;\;\;\; &\pi_2 = -\frac{3}{16} \beta_1^4 \phi^8
\end{align*}
On $\wt{\WW}$, we have $\b_1,\b_2\neq 0$, so one can see that all such points are in the orbit of a single point under the action of $(\beta_1\phi, \beta_1 \phi^2)$. This point has all of its coordinates non-zero, so it is GIT stable.
 \ed
 
 \begin{rem}\label{C0-coord-rem}
 	One can check with a computer that if one plugs in $\a_{12} = -1$, $\a_{21}=1$, $\beta_{12} = \beta_{21}= \frac32$, $\gamma_{12}= \gamma_{21}= \vareps_{21}= \frac 54$, $\vareps_{12} = -\frac54$, $\pi_i = -\frac{3}{16}$ to \eqref{g-2-n-2-curve-eq}, one obtains the union of two osculating $\P^1$'s.
	In Theorem \ref{thm:Z} below we will prove this in a different way.
\end{rem}

\subsection{Singularities of $\ov{\MM}_{2,2}(\ZZ)$}

Below we consider Smyth's $\ZZ$-stable curves, where $\ZZ$ is the extremal assigment of the unmarked components (see \cite{Smyth-modular}). The $\ZZ$-stable curves are pointed curves 
$C$ for which there exists a Deligne-Mumford stable curve $C'$ 
and a map of pointed curves $C'\to C$, contracting precisely the unmarked components
of $C'$, in a certain controlled way (see \cite[Def.\ 1.8]{Smyth-modular} for details).
We are interested in the stack $\ov{\MM}_{2,2}(\ZZ)$ of $\ZZ$-stable curves of genus $2$ with $2$ marked points.
Note that $\ov{\MM}_{2,2}(\ZZ)$ is an irreducible proper DM-stack of dimension $5$.


\begin{lem}\label{Z-stable-red-lem}
Let $(C,p_1,p_2)$ be a reducible $\ZZ$-stable point, $C_i$ the irreducible component containing $p_i$ for $i=1,2$,
$\xi=C_1\cap C_2$. Assume that $\xi$ is supported at one point then $\ell(\xi)=1$.
\end{lem}

\Pf . Indeed, if $\ell(\xi)\ge 2$ then $g(C_1)+g(C_2)\le 1$, so there exists $i$ such that $C_i\simeq\P^1$.
Since $C_i$ has only two special points, $p_i$ and the support of $\xi$, our curve cannot be $\ZZ$-stable.
\ed

\begin{prop} \label{lem:sing-codim}
(i) The only singular points of $\ov{\MM}_{2,2}(\ZZ)$ correspond to the curves $C=C_1\cup C_2$, where $p_1\in C_1$, $p_2\in C_2$,
$C_1$ and $C_2$ are both nodal curves of arithmetic genus $1$, joined transversally at their node (thus forming the singularity equivalent to the coordinate cross in the $4$-space).

\noindent (ii)	
Away from a closed subset of codimension $\ge 2$, every curve in $\ov{\MM}_{2,2}(\ZZ)$ is either smooth, or an irreducible nodal curve with normalization of genus 1, or reducible with two components of genus 1 intersecting in a node. 
\end{prop}

\Pf. (i) We use the classification of genus $2$ curves with at most two irreducible components, together with the well known fact that
the smoothness of the moduli stack near some $(C,p_1,p_2)$ is determined by the smoothness of the versal deformation spaces of the singularities of $C$.

First, let us consider the case when $C$ is irreducible. Looking at the list of possible curves (see \cite[Sec.\ 2.3]{P-contr}), we see that 
the only singularities occurring are: the node; the cusp of genus $1$, $\C[t^2,t^3]$; the tacnode; the transversal union of a cusp and a line; the coordinate cross in $3$-space;
the two cusps of genus $2$, $\C[t^3,t^4,t^5]$ and $\C[t^2,t^5]$. For all of them except perhaps for the last two the smoothness of the versal deformation spaes is well known
(for the transversal union of a cusp and a line it follows from the smoothness of the moduli space $\wt{\UU}^{ns}_{1,2}(1,0)$; see \cite[Prop.\ 3.1.1]{P-krich}).
For the two cusps of genus $2$ the required smoothness follows from the smoothness of the moduli space $\wt{\UU}^{ns}_{2,1}(2)$ proved in \cite[Prop.\ 2.1.1]{P-contr}.

Next, let us consider reducible curves $C=C_1\cup C_2$, where $p_1\in C_1$, $p_2\in C_2$. Let $\xi=C_1\cap C_2$.
As in the proof of Theorem \ref{thm:geo-stab1} we should consider the following cases.

\noindent
{\bf Case $g(C_1)=g(C_2)=1$, $\ell(\xi)=1$}. If the intersection point $q$ is smooth on one of the components then the arising singularity is a transversal union of a genus $1$
singularity with a line, so it is smooth. The case when $q$ is a cusp on one of the components does not occur since the corresponding curve $(C,p_1,p_2)$ is not $\ZZ$-stable.
There remains the case when $q$ is a node on both curves, which indeed gives a singular point of $\ov{\MM}_{2,2}(\ZZ)$.

\noindent
{\bf Case $g(C_1)=1$, $g(C_2)=0$, $\ell(\xi)=2$}. Assume first that $\xi$ is supported at two distinct points. Then these points are smooth on $C_2=\P^1$, so we again get
a transversal union of genus $1$ singularities with a line. On the other hand, the case when $\xi$ is supported at one point does not occur since such a curve is not $\ZZ$-stable.

\noindent
{\bf Case $g(C_1)=g(C_2)=0$, $\ell(\xi)=3$}. As before, by $\ZZ$-stability, it is enough to consider the case when $\xi$ is supported at more than $1$ point. If it is supported
at $3$ points then we just get $3$ nodes. If $\xi$ is supported at $2$ points then we get one node and one tacnode, so all these singularities are smooth.

\noindent
(ii) We have to go through all the strata of curves in $\ov{\MM}_{2,2}(\ZZ)$, other than the ones listed, and check that they have dimension $\le 3$. For irreducible curves 
this follows immediately from the list in \cite[Sec.\ 2.3]{P-contr}. More precisely, it is enough to check that the underlying curve depends on at most $1$ parameter.
The two cases when the dependence on $1$ parameter occurs are an elliptic curve with a cusp and $\P^1$ with two pairs of points glued nodally. All other curves do not vary.

Now let us consider the stata with $C=C_1\cup C_2$. If $g(C_1)=g(C_2)=1$, $\ell(\xi)=1$, and say, $C_1$ is singular, then the data $(C_1,p_1,q)$ is $1$-dimensional,
while the data $(C_2,p_2,q)$ is $2$-dimensional, so we get a $3$-dimensional stratum. If $g(C_1)=1$, $g(C_2)=0$, $\xi=q_1\cup q_2$, then the
data $(C_1,p_1,q_1,q_2)$ is $3$-dimensional, while the data $(C_2,p_1,q_1,q_2)$ is $0$-dimensional, so we get a $3$-dimensional stratum.
Finally, it is easy to see that the case $g(C_1)=g(C_2)=0$, $\ell(\xi)=3$ gives a $0$-dimensional stratum. 
\ed

\begin{cor}\label{Z-sm-cor} 
The stack $\ov{\MM}_{2,2}(\ZZ)$ is the union of two open substacks:
$$\ov{\MM}_{2,2}(\ZZ)=\ov{\MM}_{2,2}(\ZZ)^{sm}\cup (\ov{\MM}_{2,2}(\ZZ)\setminus (h^1(p_1+p_2)\neq 0)),$$
where $\ov{\MM}_{2,2}(\ZZ)^{sm}$ is the smooth locus in $\ov{\MM}_{2,2}(\ZZ)$, and 
the second open substack is the complement to the closed substack given by $h^1(p_1+p_2)\neq 0$.
\end{cor}

\Pf . By Proposition \ref{lem:sing-codim}, singular points of $\ov{\MM}_{2,2}(\ZZ)$ correspond to $C=C_1\cup C_2$, where $p_i\in C_i$, $g(C_1)=g(C_2)$,
with $C_1$ and $C_2$ nodal and glued transversally at the node. But for such a curve $H^0(C,\OO(p_1+p_2))=\C$, so it is in the second open substack above.
\ed

\subsection{$\ZZ$-stability versus (i)-stability}

\begin{lem}\label{GIT-vs-Z-stab-lem}
Every point in $\UU_{2,2}^{ns}$, with the exception of a single point
$[C^0]$ which is the union of two osculating $\P^1$'s, is $\ZZ$-stable.
\end{lem}

\Pf . This is easily checked by going through the list of (i)-stable curves in Theorem \ref{thm:geo-stab1}.
\ed
 



\begin{lem}\label{P1-lem} 
	(i) Let $C$ be an irreducible reduced projective curve, $p\in C$ a smooth point. If $h^0(p)>1$ then $C\simeq \P^1$.
	
	\noindent
	(ii) If $(C,p_1,p_2)$ is in $\ov{\MM}_{2,2}(\ZZ)$ then $h^0(p_1)=h^0(p_2)=1$.
\end{lem}

\Pf. (i) Let $\pi:\wt{C}\to C$ be the normalization map, and let $\wt{p}\in \wt{C}$ be the unique point over $p$.
Then $\pi^*\OO_C(p)\simeq \pi^*\OO(\wt{p})$, so $h^0(\wt{C},\wt{p})>1$. This implies that $\wt{C}\simeq \P^1$.
On the other hand, $\OO_C(p)$ is generated by global sections so it gives a degree $1$ map $C\to \P^1$
which is inverse to $\pi$. Hence, $C\simeq \wt{C}\simeq\P^1$.

\noindent
(ii) Assume that $h^0(p_1)>1$.
By (i), we see that $C$ has to be reducible. Let $C_i$ be the irreducible component containing $p_i$, for $i=1,2$,
and let $\xi=C_1\cap C_2$. By (i) we see that $C_1\simeq \P^1$ and $\xi$ is supported at one point. Thus,
$(C,p_1,p_2)$ cannot be $\ZZ$-stable. 
\ed

\begin{prop} Let $C$ be a reduced projective curve of arithmetic genus $2$, $p_1\neq p_2$ a  pair of smooth points
	such that $\OO_C(p_1+p_2)$ is ample. Assume that $h^1(2p_1+p_2)\neq 0$. Then $C$ is the union of two
	irreducible components $C_1$ and $C_2$ joined transversally at a single point, where $p_1\in C_1$, $p_2\in C_2$,
	and one of the components is isomorphic to $\P^1$.
\end{prop}

\Pf . The assumption implies that $h^1(p_1+p_2)\neq 0$. Hence, there exists a nonconstant rational function 
$f\in H^0(C,\OO(p_1+p_2))$. 

{\bf Step 1}. We claim that $C$ is reducible.
Indeed, assume that $C$ is irreducible. Then by Lemma \ref{P1-lem}, we know that $h^0(p_1)=h^0(p_2)=1$. 
Hence, $f$ has poles of order exactly $1$
at $p_1$ and $p_2$. Also, we have $h^1(2p_1)\neq 0$. Hence,
there exists a rational function $f_1\in H^0(C,\OO(2p_1))$ that has a pole of order $2$ at $p_1$.
Then the functions $1, f_1, f$ form a basis of $H^0(C,\OO(2p_1+p_2))$. Now we observe that $f_1\cdot f$ has a pole
of order $3$ at $p_1$. Hence, we get $h^0(3p_1+p_2)=4$, or equivalently $h^1(3p_1+p_2)=1$. But this implies that
$h^1(3p_1)\neq 0$, which is impossible by \cite[Lem.\ 2.4.4]{P-contr}.

{\bf Step 2}. Let $C=C_1\cup C_2$, where $C_i$ is irreducible and $p_i\in C_i$.
Assume that $f$ has a pole of order $1$ at $p_1$. We claim that in this case $C_1\simeq \P^1$ and the subscheme
$\xi:=C_1\cap C_2\sub C$ has length $1$. Indeed, $f|_{C_1}$ is a nonconstant function in $\OO_{C_1}(p_1)$,
so by Lemma \ref{P1-lem}, $C_1\simeq \P^1$. 
If $h^0(C_2,\OO(p_2))=1$ then $f|_{C_1}$ should have a constant restriction to $\xi\sub C_1$, which is possible only if
$\xi$ has length $1$. On the other hand, if $h^0(C_2,\OO(p_2))>1$ then by Lemma \ref{P1-lem}, $C_2\simeq \P^1$.
In this case, since $C$ has genus $2$, $\xi$ should have length $3$. Hence, the restriction map
$H^0(C_1,\OO(2p_1))\to H^0(\xi,\OO)$ is an isomorphism. But then the exact sequence
$$0\to H^0(C,\OO(2p_1))\to H^0(C_1,\OO(2p_1))\oplus H^0(C_2,\OO)\to H^0(\xi,\OO)$$
implies that $h^0(C,2p_1)=1$, which is a contradiction.

{\bf Step 3}. There remains a case when $f$ is regular at $p_1$ and has a pole of order $1$ at $p_2$.
Then by Lemma \ref{P1-lem}, $C_2\simeq \P^1$. Also, $f|_{C_2}$ should have a constant restriction to $\xi\sub C_2$,
which implies that $\xi$ has length $1$.
\ed

\begin{cor}\label{cor:Z-stable-h1}
	If $(C,p_1,p_2)$ is in $\ov{\MM}_{2,2}(\ZZ)$ then $h^1(2p_1+p_2)=h^1(p_1+2p_2)=0$.
\end{cor}

\begin{lem} If $(C,p_1,p_2)$ is in $\ov{\MM}_{2,2}(\ZZ)$ and $h^1(p_1+p_2)=0$ then $(C,p_1,p_2)$ is (i)-stable. \label{lem:Z-i-stable}
\end{lem}

\Pf . Assume that $(C,p_1,p_2)$ is not (i)-stable. Then $C$ is reducible: $C=C_1\cap C_2$ with $p_i\in C_i$.
Also, the subscheme $\xi=C_1\cap C_2$ has length $\le 2$.


Case 1. $\ell(\xi)=1$. Then the point $\xi$ is non-nodal on either $C_1$ and $C_2$. Without loss of generality, say it is $C_1$. This implies that $C_1$ has arithmetic genus $1$ and $\xi$ is cusp. However, in this case $C$ is not $\ZZ$-stable
(since the normalization of $C_1$ is $\P^1$ with only two distinguished points).

Case 2. $\ell(\xi)=2$. Then $\xi$ is supported at one point and one of the curves $C_1$, $C_2$ is isomorphic to $\P^1$,
hence $(C,p_1,p_2)$ is not $\ZZ$-stable.
\ed

\begin{lem} \label{lem:surj}
	If $(C, p_1, p_2)$ is in $\ov{\MM}_{2,2}(\ZZ)$ then there exists a smooth point $p_3$ on $C$ so that $h^1(p_1 + p_3)=0$.  
\end{lem}
\Pf.
If the component of $C$ containing $p_1$ has genus greater than 0, then a generic point on that component will do for $p_3$.

So we can now assume $C$ is reducible. Let $C_i$ be the component with $p_i$ (for $i=1,2$). We may assume that $C_1$ is genus 0, and that the intersection $\xi = C_1 \cap C_2$ has at least 2 points (otherwise it would not be $\ZZ$-stable).

If the genus of $C_2$ is 1, then a generic point on $C_2$ will do: we will have only constant sections of $\OO(p_1 + p_3)|_{C_2}$ and two sections of $\OO(p_1 + p_3)|_{C_1}$ which get cut down to one since they are required to be constant on $\xi$. So $h^0(p_1 + p_3)=1$, which is what we want.

The last case is when $C_1$ and $C_2$ are both rational and and $\ell(\xi) = 3$: then a generic point on $C_1$ will do: there will be three sections of $\OO(p_1 + p_3)|_{C_1}$ which get cut down to one since they are required to be constant on $\xi$.
\ed



\begin{thm} \label{thm:Z}
Let $\ov{\WW} \subset \ov{\MM}_{2,2}(\ZZ)$ be the closure of the locus $h^1(p_1 + p_2) > 0$ in $\MM_{2,2}$. Then $\ov{\WW}$ coincides with the locus $h^1(p_1 + p_2) > 0$ in $\ov{\MM}_{2,2}(\ZZ)$.

There is a regular map
\[
	\phi_3: \ov{\MM}_{2,2}(\ZZ) \rightarrow \ov{\UU}^{ns}_{2,2}
\]
such that $\phi_3(\ov{\WW})$ is a single point  and is an isomorphism elsewhere. 

One has $\phi_3(\ov{\WW})=[C^0]$, where $C^0$ is the union of two osculating $\P^1$'s.
\end{thm}

\Pf.
First, we observe that since $\ov{\MM}_{2,2}(\ZZ)$ is smooth and the locus $h^1(p_1+p_2)\neq 0$ is a Cartier divisor
(it can be given as a degeneration locus of a morphism of bundles of the same rank), 
every irreducible component of this locus has codimension 1.
By Lemma \ref{lem:sing-codim}, to see that this locus coincides with $\ov{\WW}$ it is enough to see that the locus of $(C,p_1,p_2)$ such that $C$ is either nodal with normalization of genus 1 or reducible with two components on genus 1 joined at a node, and $h^1(p_1 +p_2) >0$, is dimension $3$ (and hence codimension $2$). 

In the irreducible case, where a genus one curve is glued at points $q_1$ and $q_2$, we see that $h^0(p_1 + p_2) = 2$ implies that there is a rational function with poles at $p_1$ and $p_2$ and vanishing at $q_1$ and $q_2$. Hence $p_1 + p_2 \sim q_1 + q_2$. There are two dimensions worth of choices for $(C, p_1, p_2)$, but then a choice of $q_1$ determines $q_2$ (as the other point in the fiber of the 2:1 map to ${\mathbb P}^1$ determined by $p_1$ and $p_2$). Hence we get dimension 3 as desired.

In the reducible case, we get $h^0(p_1 + p_2)=1$ always.

Let 
\[
V^{\ZZ} \subset \wt{\UU}_{2,3}(1,0,1)
\] 
be the open substack consisting of curves $(C,p_1, p_2, p_3)$ such that $(C, p_1, p_2)$ is $\ZZ$-stable. By Corollary \ref{cor:Z-stable-h1} and Lemma \ref{lem:Z-i-stable}, we see that in fact $V^{\ZZ} \subset {\mathcal Y} \setminus Z$. We also see that the projection $V^\ZZ \rightarrow \ov{\MM}_{2,2}(\ZZ)$ is surjective by Lemma \ref{lem:surj}.

By composing  with the map $\wt{\ff}'$ (Proposition \ref{prop:forgetful}), we obtain a map 
$\phi_3': V^\ZZ \rightarrow \ov{\UU}^{ns}_{2,2}$. 

We now need only show that this map factors through the projection $V^\ZZ \rightarrow \ov{\MM}_{2,2}(\ZZ)$.
The argument is a slight modification of the one in \cite{P-contr}. The difference is that in our case $\ov{\MM}_{2,2}(\ZZ)$ is
not smooth. However, we can use the decomposition of $\ov{\MM}_{2,2}(\ZZ)$ into a smooth locus and the complement to $\ov{\WW}$ (see Corollary \ref{Z-sm-cor}). We know that the restriction of $\phi'_3$ to the preimage of
$\ov{\MM}_{2,2}(\ZZ)\setminus\ov{\WW}$ factors through $\ov{\MM}_{2,2}(\ZZ)\setminus \ov{\WW}$. It remains to prove the same
over the preimage of the smooth locus. But then we can use the same argument as in \cite{P-contr}.

Finally, by Lemma \ref{GIT-vs-Z-stab-lem}, we have a natural regular map 
$$\ov{\UU}^{ns}_{2,2}\setminus[C^0]\to \ov{\MM}_{2,2}(\ZZ)\setminus\ov{\WW},$$
where $C^0$ is the union of two osculating $\P^1$'s.
This immediately implies that $\phi_3(\ov{\WW})=[C^0]$.
\ed

\subsection{Projectivity of $\ov{M}_{2,2}(\ZZ)$ and the Weierstrass divisor in $\ov{\MM}_{2,2}(\ZZ)$}

\begin{lem}\label{W-h1-lem}
For every $(C,p_1,p_2)\in \ov{\WW}\sub\ov{\MM}_{2,2}(\ZZ)$ one has $h^1(2p_1)=0$.
\end{lem}

\Pf . Indeed, otherwise, we have nonconstant functions $f\in H^0(C,\OO(2p_1))$ and $g\in H^0(C,\OO(p_1+p_2))$.
Furthermore, since $h^0(p_1)=1$ (see Lemma \ref{P1-lem}), $g$ has a pole of order $1$ at $p_2$. This implies
that $1$, $f$ and $g$ are linearly independent in $H^0(C,\OO(2p_1+p_2))$. Hence, $h^0(2p_1+p_2)\ge 3$, so
$h^1(2p_1+p_2)\neq 0$ in contradiction with Corollary \ref{cor:Z-stable-h1}.
\ed 

Next, we are going to prove that the coarse moduli space $\ov{M}_{2,2}(\ZZ)$ is projective.
First, we consider only
curves with $h^1(2p_1)=0$ and show that the corresponding moduli of $\ZZ$-stable curves is quasiprojective.

\begin{lem}\label{Z-st-2-0-lem} 
Let us work over $\Spec(\Q)$. Then every $(C,p_1,p_2,v_1,v_2)\in \wt{\UU}^{ns}_{2,2}(2,0)$ which is $\ZZ$-stable is 
also GIT-stable. Here
we consider the GIT-stability with respect to the $\G_m^2$-action on
$\wt{\UU}^{ns}_{2,2}(2,0)$, associated with the character $(\la_1,\la_2)\mapsto \la_1\la_2$ of $\G_m^2$.
\end{lem}

\Pf . The proof is similar to that of \cite[Thm.\ 2.4.1]{P-krich}. We use standard coordinates $\a_{ij}[p,q]$
on $\wt{\UU}^{ns}_{2,2}(2,0)$.
Looking at the $\G_m^2$-weights, we see that it is enough to check that for a $\ZZ$-stable curve one has
$(\a_{21}[-1,-1],\a_{21}[-1,-2])\neq (0,0)$ and
one of the coordinates $\a_{11}[p,q]$ with $p\le -3$, $q\ge 1$, is nonzero.
If $\a_{21}[-1,-1]=\a_{21}[-1,-2]=0$ then $h^0(p_2)>1$, which is impossible by Lemma \ref{P1-lem}(ii).
On the other hand, the coordinates $\a_{11}[p,q]$ correspond to the forgetful map
$$\wt{\UU}^{ns}_{2,2}(2,0)\to \wt{\UU}^{ns}_{2,1}(2)$$
sending $(C,p_1,p_2)$ to $(\ov{C},p_1)$ with $\ov{C}=\ov{C}=\Proj(\bigoplus_n H^0(C,\OO(np_1)))$.
If all of these coordinates vanish then $\ov{C}=C^{\cusp}(2)$ is a cuspidal curve of genus $2$. But this
is possible only if $C$ is reducible and $\xi=C_1\cap C_2$ is supported at one point. By Lemma \ref{Z-stable-red-lem},
we have $\ell(\xi)=1$, so $\ov{C}=C_1$. But this implies that $C_2\simeq\P^1$, so the curve cannot be $\ZZ$-stable.
\ed

\begin{thm}\label{projective-thm} 
Let us work over $\Spec(\Q)$. Then the coarse moduli space $\ov{M}_{2,2}(\ZZ)$ is projective.
\end{thm}

\Pf . We have a proper morphism 
$$\phi_3:\ov{M}_{2,2}(\ZZ)\to \ov{U}^{ns}_{2,2},$$
which blows the Weierstrass divisor $\ov{W}$ to a point, and is an isomorphism elsewhere.
Since $\ov{U}^{ns}_{2,2}$ is also projective, it is enough to prove that the morphism $\phi_3$ is projective. 
We know that $\phi_3$ restricts to an isomorphism 
$$\ov{M}_{2,2}(\ZZ)\setminus\ov{W}\rTo{\sim} \ov{U}^{ns}_{2,2}\setminus [C^0],$$
where $(C^0,p_1,p_2)$ is the special curve which is the union of two osculating $\P^1$'s.
Thus, it is enough to prove that the morphism
$$\phi_3^{-1}(V)\to V$$
is projective, where $V$ is an open neighborhood of $[C^0]$.
We will take $V=\ov{U}^{ns}_{2,2}((1,1),(2,0))$,
i.e., the open subset of curves with $h^1(2p_1)=0$. 
Note that $V$ is given by the inequality $\a_{12}\neq 0$, so it contains $[C^0]$ (see Rem.\ \ref{C0-coord-rem}).
Also, by Lemma \ref{W-h1-lem}, $\ov{W}$ is contained in the locus $h^1(2p_1)=0$.
Thus, $\phi_3^{-1}(V)$ is precisely the open locus in $\ov{M}_{2,2}(\ZZ)$ where $h^1(2p_1)=0$. 
By Lemma \ref{Z-st-2-0-lem}, $\phi_3^{-1}(V)$ is open
in the GIT quotient of $\wt{\UU}^{ns}_{2,2}(2,0)$ by $\G_m^2$, which is a projective scheme.
Since the morphism $\phi_3^{-1}(V)\to V$ is proper we conclude that it is projective.
\ed

Finally, we will show that the Weierstrass divisor $\ov{\WW}\sub\ov{\MM}_{2,2}(\ZZ)$ can be identified with
a weighted projective stack.
Let $\ov{\UU}^{ns}_{2,1}(2)$ denote the GIT quotient of $\wt{\UU}^{ns}_{2,1}(2)$ by $\G_m$.
As was shown in \cite[Prop.\ 2.1.1]{P-contr}, there is a natural isomorphism with the weighted projective stack,
$$\ov{\UU}^{ns}_{2,1}(2)\simeq \P(2,3,4,5,6).$$

\begin{prop}\label{Weier-prop}
The forgetful map 
$$(C,p_1,p_2)\mapsto (\ov{C},p_1), \ \text{ where } \ov{C}=\Proj(\bigoplus_n H^0(C,\OO(np_1))),$$
gives an isomorphism
$$\ov{\WW}\rTo{\sim} \ov{\UU}^{ns}_{2,1}(2).$$
\end{prop}

\Pf . As before, we view $\ov{\WW}$ as a substack of $\UU^{ns}_{2,2}(2,0)$, where $h^0(p_1+p_2)\ge 2$.
Over $\wt{\UU}^{ns}_{2,2}(2,0)$ we have elements $f\in H^0(C,\OO(3p_1))$, $h\in H^0(C,\OO(4p_1))$ and
$k\in H^0(C,\OO(5p_1))$, normalized by $f\equiv 1/t_1^3+\ldots$, $h\equiv 1/t_1^4+\ldots$, $k\equiv 1/t_1^5+\ldots$,
where $t_i$ is a parameter at $p_i$ compatible with a choice of a nonzero tangent vector at $p_i$.
Furthermore, there is a unique choice of $f,h,k$, so that one has equations 
\begin{equation}\label{g-2-n-1-eq}
\begin{array}{l}
h^2=fk+q_1h+2q_1^2+f(q_{2,0}+q_{2,1}f),\\
hk=f(q_{3,0}+q_{3,1}f+f^2)-q_1k+(q_{2,0}+q_{2,1}f)h+q_1(q_{2,0}+q_{2,1}f),\\
k^2=(q_{3,0}+q_{3,1}f+f^2)h+(q_{2,0}+q_{2,1}f)^2-2q_1(q_{3,0}+q_{3,1}f+f^2),
\end{array}
\end{equation}
for some constants $q_1, q_{2,0}, q_{2,1}, q_{3,0}, q_{3,1}$ which are coordinates on
$\wt{\UU}^{ns}_{2,1}(2)$ (see \cite[Prop.\ 2.1.1]{P-contr}).
In addition, we have a function $g\in H^0(C,\OO(p_1+p_2))$, normalized by $g\equiv 1/t_1+\ldots$.
Note that $g$ has a pole of order $1$ at $p_2$, so the ring $\OO(C\setminus\{p_1,p_2\})$ has
a basis
$$f^n, f^nh, f^nk, g^{1+n}, \ \ n\ge 0.$$
Thus, we should have relations of the form
\begin{align*}
&fg=\a g+h+a(f), \\
&hg=\b g+k+dh+b(f), \\
&kg=\ga g+e_1k+e_2h+c(f),
\end{align*}
for some constants $\a,\b,\ga,e_1,e_2$ and polynomials in $f$, $a(f),b(f),g(f)$ with $\deg(a)\le 1$, $\deg(b)\le 1$,
$\deg(c)=2$ and $c$ monic. Note that $g$ is defined up to adding a constant, and we can fix this ambiguity by requiring
that $d=0$.
Applying the Buchberger's algorithm, we find 
$$\a=0, \ \ \b=2q_1,\ \ \ga=q_{2,0}, \ \ e_1=0, \ \ e_2=q_{2,1},$$
$$a(f)=q_1, \ \ b(f)=q_{2,0}+q_{2,1}f, \ \ c(f)=q_1q_{2,1}+q_{3,0}+q_{3,1}f+f^2.$$ 

Note also that for $(C,p_1,p_2)$ in $\ov{\WW}$ we cannot have $\ov{C}=C^{\cusp}(2)$ 
(see the proof of Lemma \ref{Z-st-2-0-lem}).
Thus, the forgetful map induces a well defined morphism 
$$\ov{\WW}\to \ov{\UU}^{ns}_{2,1}(2)=(\wt{\UU}^{ns}_{2,1}(2)\setminus [C^{\cusp}(2)])/\G_m^2.$$
Furthermore, the above calculations show that it is a closed embedding.
To see that its image is dense, we note that the image contains all $(C,p_1)$ with $C$ smooth (and $p_1$ not a Weierstrass
point). Indeed, we can take $p_2=\tau(p_1)$, where $\tau$ is the hyperelliptic involution on $C$.
\ed

\begin{rem} The above proposition implies that there is an involution $\si$ on $\ov{\UU}^{ns}_{2,1}(2)$ that corresponds
to the natural involution $(C,p_1,p_2)\mapsto (C,p_2,p_1)$ of $\ov{\WW}$.
On the locus of smooth curves we have $\si(C,p_1)=(C,\tau(p_1))$, where $\tau$ is the hyperelliptic involution of $C$.
\end{rem}

\end{document}